\documentclass[11pt]{article}
\usepackage[letterpaper]{geometry}
\usepackage{fullpage}
\usepackage{epsfig}
\usepackage{amsmath} \usepackage{amsfonts} \usepackage{amssymb}
\newcommand{\REMOVED}[1]{}

\newtheorem{theorem}{Theorem}

\newtheorem{lemma}{Lemma}

\newenvironment{proof}{\noindent {\em Proof}.\ }{\proofbox\par\smallskip\par}
\newcommand{\halmos}{\rule{1ex}{1.4ex}}
\newcommand{\proofbox}{\hspace*{\fill}\mbox{$\halmos$}}

\begin{document}

\title{{\bf Introducing Ramanujan's Class Polynomials \\ in the Generation of Prime Order Elliptic Curves}}

\author{Elisavet Konstantinou
\thanks{Dept of Information and Communication Systems Engineering, Univ. of the Aegean, 83200, Samos, Greece.
}
\and
Aristides Kontogeorgis \thanks{Max-Planck-Institut f\"ur Mathematik Vivatsgasse 7
 53111 Bonn,  and  Dept of Mathematics, Univ. of the Aegean, 83200, Samos, Greece.
\newline
Emails: 
%{\tt \{ekonstantinou,kontogar\}@aegean.gr}}
{\tt ekonstantinou@aegean.gr, kontogeo@mpim-bonn.mpg.de}}}
% 
% \author{Elisavet Konstantinou\thanks{Dept of Information and Communication Systems Engineering, Univ. of the Aegean, 83200, Samos, Greece.
% }
% \and
% Aristides Kontogeorgis\thanks{Dept of Mathematics, Univ. of the Aegean, 83200, Samos, Greece.
% \newline Emails: {\tt \{ekonstantinou,kontogar\}@aegean.gr}}
% }

\maketitle

\newcounter{vcounter}
\setcounter{vcounter}{1}

\begin{abstract}
Complex Multiplication (CM) method is a frequently used method for the
generation of prime order elliptic curves (ECs)
over a prime field $\mathbb{F}_p$.
The most demanding and complex step of this method is the computation of the roots
of a special type of class polynomials, called Hilbert polynomials.
These polynonials are uniquely determined by the CM discriminant $D$.
The disadvantage of these polynomials is that they have huge coefficients and thus they need
high precision arithmetic for their construction. Alternatively, Weber polynomials can be used
in the CM method. These polynomials have much smaller coefficients and their roots can be easily
transformed to the roots of the corresponding Hilbert polynomials. However, in the case
of prime order elliptic curves, the degree of Weber polynomials is three times larger
than the degree of the corresponding Hilbert polynomials and for this reason the calculation of
their roots involves computations in
the extension field $\mathbb{F}_{p^3}$.
Recently, two other classes of polynomials, denoted by $M_{D,l}(x)$ and $M_{D,p_1,p_2}(x)$ respectively, 
were introduced which can also
be used in the generation of prime order elliptic curves.
The advantage of these polynomials is that their degree is equal to the degree of the Hilbert polynomials
and thus computations over the extension field can be avoided. 

In this paper, we propose the use of a new class of polynomials. We will call them Ramanujan polynomials
named after Srinivasa Ramanujan who was the first to compute them for few values of $D$.
We explicitly describe the algorithm
for the construction of the new polynomials, show that their degree is equal to the degree of the
corresponding Hilbert polynomials and give the necessary transformation of their roots (to the roots of the corresponding Hilbert polynomials).
Moreover, we compare (theoretically and experimentally) the efficiency of
using this new class
against the use of the aforementioned Weber, $M_{D,l}(x)$ and $M_{D,p_1,p_2}(x)$ polynomials and show that 
they clearly outweigh all of them
in the generation of prime order elliptic curves.

\end{abstract}


\section{Introduction}

%Elliptic curve cryptography emerges as a fundamental and efficient
%technologic alternative to conventional public key cryptosystems for the support of {\bf auto na to bgaloume! PKI} services.
The generation of cryptographically secure elliptic curves over prime fields is one of the most fundamental
and complex problems in elliptic curve cryptography. An elliptic curve (EC) is cryptographically
secure if its use in a cryptosystem guarantees robustness against all (currently) known attacks
(e.g. \cite{FR94,MOV93,PH78,SA98}). All these attacks can be avoided if the order of the EC possesses
certain properties.
An equally important alternative to cryptographic robustness (see e.g.,
\cite{SSK01}) requires that the order of the generated EC is a
prime number. Moreover, in certain applications it is necessary that the order of the EC is
prime \cite{BLS01}.

The most commonly used methods for the generation of ECs over prime fields
are the Complex Multiplication (CM) method
\cite{AM93,LZ94,M92} and the point counting method \cite{S95}.
In this paper we follow the first approach and study the use of
the CM method for generating ECs of prime order in $\mathbb{F}_p$.
Briefly, the CM method takes as input the order $p$ of the prime field and determines
a parameter $D$ called the CM discriminant and the order $m$ of the EC.
If the order $m$ satisfies the desired properties (e.g. is a prime number) then a class
polynomial is computed using the discriminant $D$ and the parameters of the EC are constructed
from a root modulo $p$
of this polynomial. 
%
% On the other hand, point counting algorithms first choose the parameters of the
% EC and then compute its order. If the order of the EC does not satisfy the
% required properties for cryptographic robustness, then another set of EC
% parameters is generated and the process is repeated. Although ECs with no prime order can be generated
% very fast with a point counting algorithm, the requirement of prime order severely change this
% efficiency. It has been proved \cite{GMcKee00} that the probability of finding a prime order EC from the set
% of all possible ECs over a prime field $\mathbb{F}_p$
% is, asymptotically,
% $\frac{c_p}{\log p}$, where $c_p$ is a constant depending on $p$
% and satisfying $0.44 \leq c_p \leq 0.62$. Thus, it appears that
% prime orders are not especially favored by the point counting
% approach, as also noted in~\cite{BS07,GMcKee00}.
% In contrast, CM method
% starts with a prime number (the order of the EC) and {\em then}
% constructs the parameters of the EC which means that this averse prime order possibility is
% avoided.
%
The most complex and demanding step of the CM method is the computation of the class
polynomial. The original version of the method requires  the construction of a Hilbert polynomial
whose roots can be used directly for the construction of the EC parameters. The use of any
other class polynomial necessitates the existence of a transformation that will convert
the roots of this polynomial to the roots of the corresponding Hilbert polynomial.
Class polynomials are constructed with input the discriminant $D$ and by the term ``corresponding
polynomial" we mean the polynomial that is constructed with the same $D$.
The disadvantage of Hilbert polynomials is that their coefficients grow very large as the value of
discriminant increases and thus their construction requires high precision arithmetic and can be
very inefficient even for moderate values of $D$.

To overcome these shortcomings of Hilbert polynomials, two alternatives have been proposed for the case
of prime order ECs:
either to compute them off-line in powerful machines, and store
them for subsequent use (see e.g., \cite{SSK01}), or to use alternative class
polynomials for certain values of $D$ (see e.g.,
\cite{KSZ04_icisc}) and produce the required
Hilbert roots from them. The first approach however requires
storing and handling several Hilbert polynomials
with huge coefficients and this can induce problems especially in devices with limited resources.
These problems are addressed by the second approach.

Weber or $M_{D,l}(x)$ polynomials were used in the literature for the generation
of prime order elliptic curves \cite{KSZ04_icisc}. Both types of polynomials have
much smaller coefficients than the coefficients of the corresponding Hilbert polynomials and their use
can considerably improve the efficiency of the whole CM method.
 More preciselly the logarithmic height of the coefficients of the Weber and $M_{D,l}(x)$ polynomials 
is smaller by a constant factor than the corresponding logarithmic height of the 
Hilbert polynomials. 
 Weber polynomials can be computed faster than $M_{D,l}(x)$ polynomials \cite{EM02}.
However, finding their roots requires computations in the extension field
$\mathbb{F}_{p^3}$ which makes the whole process more complicated.
The reason is that in the case of prime order ECs the discriminant $D$ must be congruent
to $3 \bmod 8$ and these values give rise to Weber polynomials with degree three times larger than the
degree of the corresponding Hilbert polynomials. Thus, one must find a root of the Weber polynomial in
the extension field $\mathbb{F}_{p^3}$ and then trasform it to a root of the Hilbert polynomial in
$\mathbb{F}_{p}$. The use of $M_{D,l}(x)$ polynomials tackles this difficulty
as their degree is equal to the degree of the Hilbert polynomials.
Furthermore, the use of Weber polynomials requires the storage of three times more coefficients and the memory
needed for this purpose can be larger than the corresponding memory required for the storage
of the $M_{D,l}(x)$ polynomials. 

In \cite{ES04} the construction of another class of polynomials was proposed. We will denote these
polynomials as $M_{D, p_1, p_2}(x)$ because their construction is based on two prime numbers $p_1$ and $p_2$.
The degree of these polynomials is equal to the degree of the Hilbert polynomials and this is
a considerable advantage against Weber polynomials. Compared to the Weber polynomials, $M_{D, p_1, p_2}(x)$
polynomials have larger coefficients for all values of $p_1$ and $p_2$, except for $p_1=3,p_2=13$ and $p_1=5,p_2 = 7$.
Moreover, the modular equations which are used for the transformation of a root 
of $M_{D, p_1, p_2}(x)$
polynomials to a root of the corresponding Hilbert polynomials have degree at least 2 in the root of Hilbert
polynomial (which makes the computations more ``heavy'') and their coefficients
are quite large (which makes their use less efficient).

In conclusion, the type of polynomial that one should use depends on 
the particular application and the value of $D$.
It is clear that finding a class of polynomials which can be constructed more efficiently
than all previously mentioned polynomials, have 
degree equal to the degree of the corresponding 
Hilbert polynomials
and
have a modular equation with degree 1 in the root of Hilbert
polynomials, will considerably improve
the performance
of the CM method
for the generation of prime order elliptic curves and 
will outweigh all previously used polynomials in every
aspect (e.g. precision requirements, storage memory, time efficiency).

%previous work
Prime order ECs defined in
various fields were also treated in \cite{B01,BS06}.
In the first, the authors used the CM method with Hilbert
polynomials~\cite{B01} for the
generation of prime order ECs over extension fields, while in the second 
the authors proposed a very efficient variant of the CM method for the construction
of prime order ECs over prime fields~\cite{BS06}.
Furthermore, a number of works appeared that compare variants of the CM method and
also present experimental results concerning the construction efficiency, such
as the work of M\"{u}ller and Paulus~\cite{MP97}, as well as the theses of
Weng~\cite{WE01} and Baier~\cite{B02}.

\paragraph{Our contribution}

Srinivasa Ramanujan (1887-1920) defined on his third notebook,
pages 392 and 393 in the pagination of \cite[vol. 2]{RamNotebooks}, the
values of five class polynomials for five different values of the discriminant $D$.
The simplicity and the small coefficients of these polynomials was remarkable.
In 1999 Bruce C. Berndt and Heng Huat Chan \cite{Berndt-Chan}
proved that if $D$ is squarefree and $D \equiv 11 \bmod {24}$ then the roots of these five  polynomials
are real units and can generate the Hilbert class field. 
Moreover, they asked for an efficient way of computing these polynomials for every discriminant $D$ (and not only for
the five values computed by Ramanujan). 
In the rest of the paper, we will call them {\em Ramanujan polynomials}.

Interpreting the theorem of Berndt and Chan (that the roots of the Ramanujan polynomials can generate
the Hilbert class field for values $D \equiv 11 \bmod {24}$), we see that
Ramanujan polynomials can be used in the CM method as the aforementioned theorem proves that there is a
transformation of their roots
to the roots of the corresponding Hilbert polynomials. In addition, as $D \equiv 11 \bmod {24} \equiv 3 \bmod 8$,
Ramanujan polynomials can be used in the generation of prime order ECs.

The contribution of this paper is threefold. Firstly, we introduce for the first time the
use of Ramanujan polynomials in the CM method by providing an efficient algorithm for their
construction for all values of the discriminant. The theory behind this
construction is based on Shimura Reciprocity Law \cite{GeeBordeaux,GeeStevenhagen} 
and all the mathematical proofs behind it are presented in \cite{KonstKonto}.
However, in the context of this paper we present a considerably simplified version of the
method described in \cite{KonstKonto} which can be equally used either by a mathematician or a practitioner
with no background in algebraic number theory and algorithmic class field theory. 
%and requires a
%strong background in algebraic number theory and algorithmic class field theory.
%In the context of this paper,
%we omit the mathematical proofs and
%describe the construction of the polynomials with an eye to applications and practitioner's needs.
%For the interested reader the mathematical proofs are presented in \cite{KonstKonto}.
%In particular, in this paper we present a considerably simplified version of the method described in
%\cite{KonstKonto}.

Secondly, we observe that Ramanujan polynomials have the same degree
with their corresponding Hilbert polynomials and hence have roots
in $\mathbb{F}_{p}$. In addition, we provide the necessary transformation of a Ramanujan
polynomial's root to a root of the corresponding Hilbert polynomial and thus give all the information
that a practioner needs in order to use the new class of polynomials in the CM method.

Finally, we perform a comparative theoretical and experimental
study regarding the efficiency of the CM method using the
aforementioned Weber, $M_{D,l}(x)$ and $M_{D,p_1,p_2}(x)$ polynomials against the new class of
polynomials.
We
show that Ramanujan polynomials are by far the best choice when CM
method is used for the generation of prime order elliptic curves because their degree is equal to
the degree of the corresponding Hilbert polynomials and their construction is more efficient
than the construction of all previously used polynomials.
We show that the logarithmic height of the coefficients of the Ramanujan polynomials is asymptotically 36 times smaller 
than the logarithmic height of the Hilbert polynomials and this allows us to 
 show that the precision requirements for the construction of Ramanujan
polynomials can be from 22\% to 66\% smaller than the precision requirements of all other
class polynomials.

In literature  the ``efficiency'' of a class invariant (a root of a class polynomial) is measured 
by the asymptotic ratio of the logarithmic height of a root of the Hilbert polynomial to 
a root of the class polynomial in question. The best known class invariant is the one used
for the construction of 
Weber polynomials with $D \not\equiv 0 \pmod 3$ and $D \equiv 3,7 \pmod 8$. 
The roots of these Weber polynomials have logarithmic height that is  asymptotically 72  times smaller  than 
the logarithmic height of the roots  of the corresponding Hilbert polynomials. 
However, in practice we are not interested in the logarithmic height 
of the roots but in the logarithmic height of the polynomials, since the latter 
measures the precision required for the construction of the polynomials. 
In this paper we will show that these two heights coincide only if the class polynomial has 
degree equal to the degree of the corresponding Hilbert polynomial. 
For the construction of prime order elliptic curves, 
Weber class polynomials have degree 3 times larger than the degree of the Hilbert polynomials. 
We will show that in this case the logarithmic height
of the Weber polynomials  is asymptotically 24=72/3 times less than 
the logarithmic height of Hilbert polynomials and not 72. Thus, even though the height of Weber polynomials' roots
is smaller than the height of the roots of Ramanujan's class polynomials, the precision requirements for the
construction of the latter are smaller.

Ramanujan polynomials can also be used in the generation of special curves, such as MNT curves \cite{MNT00,MNT01,SB04} and in the generation of ECs that do not necessarily have prime order \cite{AM93,LZ94}.
It is interesting to note here that in the latter case, as our experiments indicated, Ramanujan polynomials
%It is interesting to note here that even in this case they
outweigh Weber polynomials for all values of the discriminant $D \not\equiv 7 \bmod 8$.
Moreover, problems such as primality testing/proving \cite{AM93} and the representability
of primes by quadratic forms \cite{C89} can be considerably improved with the use of Ramanujan polynomials.
This makes our analysis for these polynomials
even more useful.

The rest of the paper is organized as follows. In
Section~\ref{pprel} we review some basic definitions and facts
about ECs and the CM method.
In Section~\ref{class_polynomials} we review properties of Hilbert, Weber, $M_{D,l}(x)$ and $M_{D,p_1,p_2}(x)$ polynomials with
$D \equiv 3 \bmod 8$ and
in Section~\ref{ramanujan} we elaborate on the construction of Ramanujan
polynomials describing in an explicit way how they can be used in the CM method. 
In Section~\ref{prec-section} we provide theoretical estimations for the precision requirements of all
previously mentioned polynomials and 
%Finally,
in Section~\ref{exper} we present our experimental results.
%Preliminary parts of this work appeared in \cite{KK06}.

\section{A Brief Overview of Elliptic Curve Theory and Complex Multiplication}
\label{pprel}

In this section we give a brief introduction to elliptic curve
theory and to the Complex Multiplication method for generating prime
order elliptic curves. Our aim is to facilitate the reading of the sections
that follow.
%For full coverage of the necessary concepts and
%terms, the interested reader may consult \cite{BSS99}.
%Also, the proofs of certain theorems require basic knowledge of
%algebraic number theory and Galois theory.
%The interested reader is referred to \cite{C89,S04,ST87} for definitions
%not given here.

\subsection{Preliminaries of Elliptic Curve Theory}
\label{prel}

An {\em elliptic curve} over a finite field $\mathbb{F}_{p}$, $p$ a prime
larger than 3, is denoted by $E(\mathbb{F}_{p})$ and it is comprised of all
the points $(x, y) \in \mathbb{F}_{p}$ (in affine coordinates) such that
\begin{equation}
y^2 = x^3 + ax + b, \label{ec}
\end{equation}
with $a, b \in \mathbb{F}_{p}$ satisfying $4a^3 + 27b^2 \neq 0$. These
points, together with a special point denoted by $\cal O$ (the
{\em point at infinity}) and a properly defined addition operation
form an Abelian group. This is the {\em Elliptic Curve group} and
the point $\cal O$ is its zero element (see \cite{ACDFLNV06,BSS99,S86}
for more details on this group).

The {\em order}, denoted by $m$, is the number of points that
belong in $E(\mathbb{F}_{p})$.
The difference between $m$ and $p$ is measured by the so-called
{\em Frobenius trace} $t=p+1-m$ for which Hasse's theorem (see e.g.,
\cite{BSS99}) states that $|t|\leq 2\sqrt{p}$, implying that
$p + 1 - 2\sqrt{p} \leq m \leq p + 1 + 2\sqrt{p}$.
%
%\begin{equation*}
%p + 1 - 2\sqrt{p} \leq m \leq p + 1 + 2\sqrt{p}.
%\label{eq:hasse}
%\end{equation*}
%
This is an important inequality that provides lower and upper bounds on the
number of points in an EC group.
The {\em order} of an element  $P\in E(\mathbb{F}_p)$ is defined as the
smallest positive integer $n$ such that $nP = \cal O$.  Langrange's
theorem implies that the order of a point $P\in E(\mathbb{F}_p)$
divides the order $m$ of the group $E(\mathbb{F}_p)$. Thus, $mP =
\cal O$ for any $P\in E(\mathbb{F}_p)$ and, consequently, the
order of a point is always less than or equal to the order of the
elliptic curve.

Among the most important quantities defined for an elliptic curve
$E(\mathbb{F}_{p})$ are the {\em curve
discriminant} $\Delta$ and the {\em $j$-invariant}. These two
quantities are given by the equations $\Delta = -16(4a^3 + 27b^2)$
and $j = -1728(4a)^3/\Delta$.
Given a $j$-invariant $j_0\in \mathbb{F}_p$
(with $j_0\neq 0,1728$) {\em two} ECs can be constructed. If $k =
j_0/(1728-j_0) \bmod p$, one of these curves is given by
Eq.~(\ref{ec}) by setting $a = 3k \bmod p$ and $b = 2k \bmod p$.
The second curve (the {\em twist} of the first) is given by the
equation $y^2 = x^3 + ac^2x + bc^3$
%
%\begin{equation*}
%    y^2 = x^3 + ac^2x + bc^3
%\label{ecp2}
%\end{equation*}
%
with $c$ any quadratic non-residue of $\mathbb{F}_{p}$.
If $m_1$ and $m_2$ denote the orders of an elliptic curve and its
twist respectively, then $m_1+m_2=2p+2$ which implies that if one
of the curves has order $p+1-t$, then its twist has order $p+1+t$,
or vice versa (see~\cite[Lemma VIII.3]{BSS99}).

\subsection{The Complex Multiplication Method}
\label{cmmethod}

As stated in the previous section, given a $j$-invariant one
may readily construct an EC. Finding a suitable $j$-invariant for
a curve that has a given order $m$ can be accomplished through the
theory of {\em Complex Multiplication} (CM) of elliptic curves
over the rationals. This method is called the {\em CM method} and
in what follows we will give a brief account of it.

By Hasse's theorem, $Z = 4p - (p+1-m)^2$ must be positive and,
thus, there is a unique factorization $Z = Dv^2$, with $D$ a
square free positive integer. Therefore
\begin{equation}
  4p = u^2 + Dv^2
\label{eq:D}
\end{equation}
for some integer $u$ that satisfies the equation
\begin{equation}
  m = p + 1 \pm u.
\label{orderm}
\end{equation}
The negative parameter $-D$ is called a {\em CM discriminant for
the prime $p$}. For convenience throughout the paper, we will use
(the positive integer) $D$ to refer to the CM discriminant.
The CM method uses $D$ to determine a $j$-invariant. This
$j$-invariant in turn, will lead to the construction of an EC of
order $p+1-u$ or $p+1+u$.

The CM method works as follows. Given a prime $p$, the
smallest $D$ is chosen for which there exists some integer $u$ for
which Eq.~(\ref{eq:D}) holds.
If neither of the possible orders $p+1-u$ and $p+1+u$ is
suitable for our purposes, the process is repeated with a new
$D$. If at least one of these orders is suitable, then
the method proceeds with the construction of the {\em Hilbert
polynomial} (uniquely defined by $D$) and the determination of
its roots modulo $p$.
Any root of the Hilbert polynomial can be used as a
$j$-invariant. From this root the corresponding
EC and its twist can be constructed as described in Section~\ref{prel}.
In order to find which one of the curves has the desired suitable
order ($m=p+1-u$ or $m=p+1+u$), Langrange's
theorem can be used as follows: we repeatedly choose points $P$ at random in
each EC until a point is found in one of the curves for which
$mP\neq{\cal O}$. This implies that the curve we seek is the other
one. Recently, different methods have been proposed for
choosing efficiently the correct elliptic curve in CM method \cite{NM05,RS07}.

The most demanding step of the CM method is the construction of the Hilbert polynomial, as it
requires high precision floating point and complex arithmetic.
As the value of the discriminant $D$ increases, the coefficients of the polynomials grow extremely
large and their computation becomes more inefficient.
In~\cite{B02,KSZ02}, a variant of the CM method was proposed to avoid this problem.
This variant starts with
a discriminant $D$ and a specific prime $p$ chosen at random, or
from a set of prescribed primes. It then computes $u$ and $v$
using Cornacchia's algorithm~\cite{C08} to solve Eq.~(\ref{eq:D}), and requires that the resulting
EC order $m$ is suitable (cf.~Section \ref{prel}).
Using this variant, the user can choose the value of the discriminant he wishes (and thus avoid very large values which was not possible in the original version of the CM method) or he can construct the Hilbert polynomials in a preprocessing phase and store them for later use.
In this way, the burden of their costly computation can be avoided during
the execution of the CM method. A similar variant was proposed in \cite{SSK01} for the construction of prime
order ECs.

We now turn to the generation of prime order ECs.
If $m$ should be a prime
number, then it is obvious that $u$ should be odd. It is also easy
to show that $D$ should be congruent to $3 \bmod 8$ and $v$ should
be odd, too.
In this paper, we follow the variant of the CM method proposed in ~\cite{B02,KSZ02} for the construction
of prime order elliptic curves. Thus,
we start with a CM discriminant $D \equiv 3 \bmod 8$ for the
computation of the Hilbert polynomial,
and then generate at random, or select from a pool
of precomputed {\em good} primes (e.g., Mersenne primes),
a prime $p$ and compute odd integers $u, v$ such that $4p =
u^{2}+Dv^{2}$.
Those odd integers $u, v$ can be computed with four different ways,
which are outlined in \cite{KSZ04_icisc}.
Once we have found primes $p$ and $m$ which satisfy Eq.~(\ref{eq:D}) and Eq.~(\ref{orderm}),
we can
%The most trivial
%way to do this, is by choosing at random odd integers $u$ and $v$ and then check
%whether $p$ and $m$ are prime using Eq.~(\ref{eq:D}) and Eq.~(\ref{orderm}).
%If no such numbers $u$ and $v$ can be found,
%then take another prime $p$ and repeat. Otherwise,
proceed with the next steps, which
are similar to those of the original CM method.

If we could find a way to compute the roots of the Hilbert polynomials directly, it is clear that it
wouldn't be necessary to construct the polynomials (since only their roots are needed in the CM method).
%It is clear that in CM method it will not be necessary to compute the Hilbert polynomials if we can find
%a way to compute their roots directly (since $j$-invariants are roots of the
%Hilbert polynomials).
Indeed, there are polynomials (known as class polynomials) \cite{EM02,ES03,KVY89,S02}
with much
smaller coefficients, which can be
constructed much more efficiently than Hilbert polynomials and their roots can be transformed to
the roots of the Hilbert polynomials. Thus, we can replace the Hilbert polynomials in the CM method
with another class of polynomials given that their roots can be transformed to the roots of the Hilbert
polynomials.
% In particular for the case of prime order ECs, the use of
% Weber and $M_{D,l}(x)$ polynomials was proposed in \cite{KSZ04_icisc}. 
In the following section we will briefly review the definition of these polynomials along with another class of polynomials defined in \cite{ES04} (denoted as $M_{D,p_1,p_2}(x)$) and show how they can be
used in the CM method, while in Section \ref{ramanujan} we will propose the use of Ramanujan class
polynomials.

\section{Class Polynomials}
\label{class_polynomials}

In this section we define Hilbert, Weber, $M_{D,l}(x)$ and $M_{D,p_1,p_2}(x)$ polynomials for discriminant
values $D \equiv 3 \bmod 8$ and briefly discuss their use in the CM method. The interested
reader is referred to \cite{ES04,KSZ04_icisc} for proofs and details not given here.

\subsection{Hilbert Polynomials}
\label{sec_hpoly}

Every CM discriminant $D$ defines a unique Hilbert polynomial,
denoted by $H_{D}(x)$. Given a positive $D$, the Hilbert
polynomial $H_{D}(x) \in \mathbb{Z}[x]$ is defined as
\begin{equation}
  H_{D}(x) = \prod_{\tau} (x-j(\tau))
\label{hx}
\end{equation}
for values of $\tau$ satisfying $\tau = (-\beta +
\sqrt{-D})/2\alpha$, for all integers $\alpha$, $\beta$, and
$\gamma$ such that (i) $\beta^2 - 4\alpha\gamma = -D$, (ii)
$|\beta| \leq \alpha \leq \sqrt{D/3}$, (iii) $\alpha \leq \gamma$,
(iv) $\gcd(\alpha, \beta, \gamma) = 1$, and (v) if $|\beta| =
\alpha$ or $\alpha = \gamma$, then $\beta \geq 0$.
The 3-tuple of integers $\left[ \alpha, \beta, \gamma \right]$
that satisfies these conditions is called a {\em primitive,
reduced quadratic form} of $-D$, with $\tau$ being a root of the
quadratic equation $\alpha z^{2}+\beta z+\gamma=0$. Clearly, the
set of primitive reduced quadratic forms of a given discriminant
is finite. The quantity $j(\tau)$ in Eq.~(\ref{hx}) is called {\em
class invariant} and is defined as follows. Let $z =
e^{2\pi\sqrt{-1}\tau}$ and $h(\tau) = \left(
\frac{\eta(2\tau)}{\eta(\tau)} \right)^{24}$, where $\eta(\tau)
= z^{1/{24}}\left( 1+\sum_{n \geq 1}
  {(-1)^n\left(z^{n(3n-1)/2}+z^{n(3n+1)/2}\right)}\right)$ is the Dedekind eta-function.
Then, $j(\tau) = \frac{(256h(\tau)+1)^3}{h(\tau)}$.
It can be shown \cite{C89} that Hilbert polynomials with degree
$h$ have $h$ roots modulo $p$ when they are used in the CM method.

\REMOVED{
Let $h$ be the number of primitive reduced quadratic forms, which
determines the {\em degree} (or {\em class number}) of
$H_D(x)$. Then, the bit precision required for the
generation of $H_D(x)$ can be estimated
(see \cite{LZ94}) by
%
%\begin{equation*}
$\mbox{H-Prec}(D) \approx \frac{\ln 10}{\ln 2} (h/4 + 5) +
\frac{\pi\sqrt{D}}{\ln 2} \sum_{\tau} \frac{1}{\alpha}
$
%\label{hplz}
%\end{equation*}
%
with the sum running over the same values of $\tau$ as the product
in Eq.~(\ref{hx}). It can be shown \cite{C89} that Hilbert polynomials with degree
$h$ have $h$ roots modulo $p$ when they are used in the CM method.
%under certain conditions stated in the following theorem.
}

\REMOVED{

\begin{theorem} \cite{KSZ04_icisc}
A Hilbert polynomial $H_D(x)$ with degree $h$ has exactly $h$ roots modulo $p$ if and only if
the equation $4p = u^2 + Dv^2$ has integer solutions and $p$ does not divide the
discriminant\footnote{For a definition of the discriminant of a
polynomial see \cite{C93}.} $\Delta(H_{D})$ of the polynomial.
\label{hilbertroots}
\end{theorem}

There are finitely many primes dividing the discriminant $\Delta(H_{D})$
of the Hilbert polynomial and infinitely many primes to choose. In
elliptic curve cryptosystems the prime $p$ is at least 160 bits.
Therefore, an arbitrary prime almost certainly does not divide the
discriminant. In addition, the above theorem indicates that from a Hilbert polynomial
with $h$ roots, one can readily construct $h$ ECs using the CM method.
}

\subsection{Weber Polynomials}

The Weber polynomial $W_{D}(x)\in\mathbb{Z}[x]$ for $D \equiv 3 \bmod 8$
is defined as
\begin{equation*}
  W_{D}(x) = \prod_{\ell} (x-g(\ell))
%\label{wx}
\end{equation*}
where $\ell = \frac{-b+\sqrt{-D}}{a}$ satisfies the equation $ay^2 + 2by + c=0$
for which $b^2 - ac= -D$ and (i) $\gcd(a,b,c)=1$, (ii) $|2b| \leq a \leq c$, and
(iii) if either $a=|2b|$ or $a=c$, then $b \geq 0$.
Let $\zeta=e^{\pi\sqrt{-1}/24}$. 
The class invariant $g(\ell)$ for $W_D(x)$ is defined by 
\begin{equation*}
g(\ell)= \left\{
               \begin{array}{rl}
               \zeta^{b(c - a - a^{2}c)} \cdot f(\ell) & \mbox{ if $2\mid\!\!\!\!/a$ and $2\mid\!\!\!\!/c$} \\
               -(-1)^{\frac{a^2-1}{8}} \cdot \zeta^{b(ac^2 - a -2c)} \cdot f_1(\ell)
                                                      & \mbox{ if $2\mid\!\!\!\!/a$ and $2\mid c$} \\
               -(-1)^{\frac{c^2-1}{8}} \cdot \zeta^{b(c - a - 5ac^2)} \cdot f_2(\ell)
                                                      & \mbox{ if $2\mid a$ and $2\mid\!\!\!\!/ c$}
               \end{array}
               \right.
%\label{gl1}
\end{equation*}
if $D \equiv 3 \bmod 8$ and $D \not\equiv 0 \bmod 3$, and

\begin{equation*}
g(\ell)= \left\{
               \begin{array}{rl}
               \frac{1}{2}\zeta^{3b(c - a - a^{2}c)} \cdot f^3(\ell) & \mbox{ if $2\mid\!\!\!\!/a$
               and $2\mid\!\!\!\!/c$} \\
               -\frac{1}{2}(-1)^{\frac{3(a^2-1)}{8}} \cdot \zeta^{3b(ac^2 - a -2c)} \cdot f_1^3(\ell)
                                                      & \mbox{ if $2\mid\!\!\!\!/a$ and $2\mid c$} \\
               -\frac{1}{2}(-1)^{\frac{3(c^2-1)}{8}} \cdot \zeta^{3b(c - a - 5ac^2)} \cdot f_2^3(\ell)
                                                      & \mbox{ if $2\mid a$ and $2\mid\!\!\!\!/ c$}
               \end{array}
               \right.
%\label{gl2}
\end{equation*}
if $D \equiv 3 \bmod 8$ and $D \equiv 0 \bmod 3$.
The functions $f()$, $f_1()$ and $f_2()$ are called Weber functions and are defined by
(see \cite{AM93,ieee}):
\begin{eqnarray*}
f(y) & = & q^{-1/48}\prod_{r=1}^{\infty}(1+q^{(r-1)/2}) ~~~~~~~
f_1(y)  = q^{-1/48}\prod_{r=1}^{\infty}(1-q^{(r-1)/2}) \\
f_2(y) & = & \sqrt{2}~~q^{1/24}\prod_{r=1}^{\infty}(1+q^{r}) ~~~~~~~
\mbox{ where } q = e^{2\pi y\sqrt{-1}}.
\end{eqnarray*}

For these cases of the discriminant ($D \equiv 3 \bmod 8$), the
Weber polynomial $W_D(x)$ has degree three times larger than the degree
of its corresponding Hilbert polynomial $H_D(x)$.
In \cite{KSZ04_icisc} it is shown that the Weber polynomial has roots in the
extension field $\mathbb{F}_{p^3}$.
Thus, in order to use Weber polynomials in the CM
method we must find at least one of their roots in the extension field $\mathbb{F}_{p^3}$.
The idea is that we replace Hilbert polynomials with
Weber polynomials and then try to compute a root of the Hilbert polynomial from
a root of its corresponding Weber polynomial.
To compute the desired Hilbert root, we proceed in three stages. First, we construct
the corresponding Weber polynomial. Second, we compute its roots in $\mathbb{F}_{p^3}$.
Finally, we transform the Weber roots to the desired Hilbert roots in $\mathbb{F}_p$ using
a modular equation $\Phi_{W}(x,j) = 0$. In particular, if $x$ is a root of Weber polynomial
and $j$ is a root of the corresponding Hilbert polynomial, then
\begin{equation}
\label{eq-phi-1}
\Phi_{W}(x,j) = (2^{12}x^{-24}-16)^3-2^{12}x^{-24}j
\end{equation}
if $D \not\equiv 0 \pmod 3$ and
\begin{equation}
\label{eq-phi-2}
\Phi_{W}(x,j) = (2^{4}x^{-8}-16)^3-2^{4}x^{-8}j
\end{equation}
if $D \equiv 0 \pmod 3$.
To compute a root of $W_D(x)$ in $\mathbb{F}_{p^3}$, we have to find an irreducible factor
(modulo $p$) of degree 3 of the polynomial. This can be achieved using
Algorithm 3.4.6 from \cite{C93}. The irreducible factor has 3 roots in $\mathbb{F}_{p^3}$
from which it suffices to choose one, in order to accomplish the third stage.
Details on the use of Weber polynomials in the construction of prime order elliptic curves can be
found in \cite{KSZ04_icisc}.

%An upper bound for the precision requirements of Weber polynomials
%for both cases of $D$ was presented in \cite{KSZ03} and is equal to
%$3h + \frac{\pi\sqrt{D}}{24\ln{2}} \sum_{\ell}\frac{1}{\alpha}$ for
%$D \not\equiv 0 \bmod 3$ and to $3h + \frac{\pi\sqrt{D}}{8\ln{2}}
%\sum_{\ell}\frac{1}{\alpha}$ for $D \equiv 0 \bmod 3$.
%The sum runs over the same values of $\ell$ as the product
%of Eq.~(\ref{wx}) and $3h$ is the degree of the Weber polynomial
%($h$ is the degree of the corresponding Hilbert polynomial).

\subsection{$M_{D,l}(x)$ Polynomials}
\label{mdl_poly}

Even though Weber polynomials have much smaller coefficients than Hilbert polynomials
and can be computed very efficiently, the fact that their degree for $D \equiv 3 \bmod 8$
is three times larger than the degree of the corresponding Hilbert polynomials can be a potential problem,
because it involves computations in extension fields.
Moreover, the computation of a cubic factor modulo $p$ in a polynomial with degree $3h$ is more time consuming
than the computation of a single root modulo $p$ of a polynomial with degree $h$.

To alleviate these problems, the use of a relatively new class of
polynomials was proposed 
%in \cite{KSZ04_icisc} 
referred as the $M_{D,l}(x)$ polynomials.
These polynomials have degree $h$ like Hilbert polynomials and thus they have roots modulo $p$.
They are constructed from a family of $\eta$-products:
$m_l(z)=\frac{\eta(z/l)}{\eta(z)}$ \cite{M00} for an integer $l \in \left\{3,5,7,13\right\}$.
The polynomials are obtained from this
family by evaluating their value at a suitably chosen system of quadratic forms. Once a polynomial is
computed, we can use a modular equation $\Phi_{l}(x,j) = 0$ (see Table~\ref{f_l}),
in order to compute a root $j$ modulo $p$ of the Hilbert polynomial
from a root $x$ modulo $p$ of the $M_{D,l}(x)$ polynomial.
\begin{table}
\begin{center}
\begin{tabular}{|p{1.5cm}|p{7.5cm}|} \hline
 $l$ & $\Phi_l(x,j)$ \\ \hline \hline
 3 & $(x+27)(x+3)^3-jx$ \\ \hline
 5 & $(x^2+10x+5)^3-jx$ \\ \hline
 7 & $(x^2+13x+49)(x^2+5x+1)^3-jx$ \\ \hline
 13 & $(x^2+5x+13)(x^4+7x^3+20x^2+19x+1)^3-jx$ \\ \hline
\end{tabular}
\end{center}
\caption{Modular functions for different values of $l$.}
\label{f_l}
\end{table}
%
 
%The interested reader is referred to
%\cite{KSZ04_icisc} for more details on the construction of the polynomials and on their use in the CM method
%for the generation of prime order elliptic curves.
%
%We finally note that the precision required for the construction of the
%$M_{D, l}(x)$ polynomials is approximately $\frac{1}{(l+1)}\frac{\pi\sqrt{D}}{\ln 2} \sum_{\tau} \frac{1}{\alpha}$
%where the sum runs over the same values of $\tau$ as the product
%in Eq.~(\ref{hx}) \cite{EM02}.

\subsection{$M_{D,p_1,p_2}(x)$ Polynomials}
\label{mdl_double_poly}

In authors of \cite{ES04} proposed the use of another class of polynomials. Like $M_{D,l}(x)$ polynomials,
these polynomials are constructed 
using a family of $\eta$-products:
$m_{p_1,p_2}(z) = \frac{\eta(z/p_1)\eta(z/p_2)}{\eta(z/(p_1p_2))\eta(z)}$.
%where $p_1, p_2$ are primes such that $24 | (p_1-1)(p_2-1)$.
We will refer to the minimal polynomials of these products 
%(powers of which generate
%the Hilbert class field and are called class invariants like $j(\tau)$)
as $M_{D, p_1,p_2}(x)$ where $D$ is the discriminant used for
their construction. 
The only restriction posed on the discriminant is that 
$
\left( 
\frac{D}{p_1}
\right)\neq -1$ and $\left( 
\frac{D}{p_2}
\right)\neq -1$ if $p_1\neq p_2$ or 
$\left( 
\frac{D}{p}
\right)\neq -1$ if $p_1=p_2=p$, where
$\left( \frac{\cdot}{\cdot}
\right)$ is the symbol of Kronecker.
The polynomials are obtained from this
family of $\eta$-products by evaluating their value at a suitably chosen system of quadratic forms. 
In particular, the polynomial $M_{D,p_1,p_2}(x)\in\mathbb{Z}[x]$ 
is defined as
\begin{equation*}
  M_{D,p_1,p_2}(x) = \prod_{\tau_Q} (x-m_{p_1,p_2}(\tau_Q))
%\label{mldx}
\end{equation*}
where $\tau_Q = \frac{-B_i+\sqrt{-D}}{2A_i}$ for all representatives $S = \left\lbrace (A_i, B_i, C_i) \right\rbrace_{1 \leq i \leq h} $ of the reduced primitive quadratic forms of a discriminant $-D$ derived from a $(p_1p_2)$-system \cite{S02}.

Once a polynomial is
computed, we can use the modular equations $\Phi_{p_1,p_2}(x,j)=0$,
in order to compute a root $j$ modulo $p$ of the Hilbert polynomial
from a root $x$ modulo $p$ of the $M_{D, p_1,p_2}(x)$ polynomial.
However, a disadvantage of $M_{D, p_1,p_2}(x)$ polynomials is that 
the corresponding modular polynomials $\Phi_{p_1,p_2}(x,j)$ have
degree at least 2 in
%the degree in
$j$ (which makes the computations more ``heavy'') and their coefficients
are quite large (which makes their
use less efficient) \footnote{For example, 
notice in \cite{ES05} the size of the smallest modular polynomial $\Phi_{5,7}(x,j)$.}.  
\REMOVED{
for $p_1, p_2 \notin \left\{2,3\right\}$ is $\Phi_{5,7}(x,j)$ and is equal to:
\begin{eqnarray*}
\Phi_{5,7}(x,j) & = & x^{48} + (-j+708)x^{47} + (35j+171402)x^{46} + (-525j+15185504)x^{45}\\
&   & + (4340j+248865015)x^{44} + (-20825j+1763984952)x^{43} \\
&   & + (52507j+6992359702)x^{42} + (-22260j+19325688804)x^{41} \\
&   & + (-243035j+42055238451)x^{40} + (596085j+70108209360)x^{39} \\
&   & + (-272090j+108345969504)x^{38} + (-671132j+121198179480)x^{37} \\
&   & + (969290j+155029457048)x^{36} + (-1612065j+97918126080)x^{35} \\
&   & + (2493785j+141722714700)x^{34} + (647290j-1509796288)x^{33} \\
&   & + (-3217739j+108236157813)x^{32} + (3033590j-93954247716)x^{31} \\
&   & + (-5781615j+91135898154)x^{30} + (1744085j-108382009680)x^{29} \\
&   & + (1645840j+66862445601)x^{28} + (-2260650j-66642524048)x^{27} \\
&   & + (6807810j+38019611082)x^{26} + (-2737140j-28638526644)x^{25} \\
&   & + (2182740j+17438539150)x^{24} + (-125335j-8820058716)x^{23} \\
&   & + (-1729889j+5404139562)x^{22} + (1024275j-1967888032)x^{21} \\
&   & + (-1121960j+1183191681)x^{20} + (395675j-370697040)x^{19} \\
&   & + (-54915j+103145994)x^{18} + (15582j-42145404)x^{17} \\
&   & + (34755j-15703947)x^{16} + (-6475j-3186512)x^{15} \\
&   & + (1120j-4585140)x^{14} + (-176j+1313040)x^{13} \\
&   & + (j^2-1486j-38632)x^{12} + (-7j+399000)x^{11} \\
&   & + (-19j+211104)x^{10} + (-9j+6771)x^{8} + (8j-6084)x^{7} \\
&   & + (7j-5258)x^{6} + (j-792)x^{5} - 105x^{4} + 16x^{3} + 42x^{2} + 12x + 1\\
\end{eqnarray*}
}
The only modular polynomials that have degree 2 in $j$ are
$\Phi_{3,13}(x,j)$ and $\Phi_{5,7}(x,j)$.
In addition, $M_{D, 3, 13}(x)$ and $M_{D, 5, 7}(x)$ polynomials
are constructed more efficiently than other polynomials of the double eta family \cite{EM02}.
Thus, we only used these polynomials in our experiments.  
%Finally, we note that the precision required for the construction of the
%$M_{D, p_1, p_2}(x)$ polynomials is approximately $\frac{(p_1-1)(p_2-1)}{12(p_1+1)(p_2+1)}\frac{\pi\sqrt{D}}{\ln 2} \sum_{\tau} \frac{1}{\alpha}$
%where the sum runs over the same values of $\tau$ as the product
%in Eq.~(\ref{hx}) \cite{EM02}. Thus, it is equal to $\frac{1}{28}\frac{\pi\sqrt{D}}{\ln 2} \sum_{\tau} \frac{1}{\alpha}$ for $M_{D, 3, 13}(x)$ polynomials and to $\frac{1}{24}\frac{\pi\sqrt{D}}{\ln 2} \sum_{\tau} \frac{1}{\alpha}$ for $M_{D, 5, 7}(x)$ polynomials. 

\section{Ramanujan Polynomials}
\label{ramanujan}

In this section, we define a new class of polynomials which can be used in the CM method for the
generation of prime order ECs. We elaborate on their construction and provide the necessary
transformations of their roots to the roots of the corresponding Hilbert polynomials.

\subsection{Construction of Polynomials}

Srinivasa Ramanujan (1887-1920) defined
on his third notebook, pages 392 and 393 in the pagination of \cite[vol. 2]{RamNotebooks}
the values
\begin{equation} \label{tndef}
t_D=\sqrt{3} q_D^{1/18} \frac{f(q_D^{1/3}) f(q_D^3)}{f^2(q_D)} \in \mathbb{R}
\end{equation}
where $f(-q)=\prod_{d=1}^\infty (1-q^d) = q^{-1/24}\eta(\tau)$,
%\begin{equation*}\label{Ramt}
%f(-q)=\prod_{d=1}^\infty (1-q^d) = q^{-1/24}\eta(\tau),
%\end{equation*}
$q= \exp(2\pi i \tau)$, $q_D=\exp(-\pi \sqrt{D})$, $\tau \in \mathbb{H}$ ($\mathbb{H}$ is the upper half plane) and $\eta(\tau)$ denotes the Dedekind eta-function.
Without any further explanation on how he found them,
 Ramanujan gave the following table of polynomials $T_D(x)$ based on $t_D$ for five values of $D$:
\[
\begin{array}{|c|c|}
\hline
 D & T_D(x)\\
\hline
11 & x-1 \\
35 & x^2+x-1\\
59 & x^3+2x-1 \\
83 & x^3+2x^2+2x-1\\
107 & x^3-2x^2+4x-1\\
\hline
\end{array}
\]
In \cite{Berndt-Chan}  Bruce C.  Berndt and Heng Huat Chan   proved  that these polynomials indeed
have roots the Ramanujan values $t_D$. The method they used could not be applied for
higher values of $D$ and they asked for an efficient way of computing the polynomials $T_D$ for every $D$.
They also proved that if $D \in \mathbb{N}$ is squarefree so that $D\equiv 11 \bmod {24}$ then
$t_D$ is a real unit generating the Hilbert class field.
This actually means that the polynomials $T_D$ can be used in the CM method because their roots can be transformed
to the roots of the corresponding Hilbert polynomials. In addition, the remarkably small coefficients
of these polynomials are a clear indication that their use in the CM method can be especially favoured.

In this paper we will elaborate on the construction of these polynomials, which we will call {\em Ramanujan polynomials}
and we will provide an efficient algorithm for their computation
for every discriminant $D \equiv 11 \bmod {24}$.
The theory behind this
construction is based on Shimura Reciprocity Law \cite{GeeBordeaux,GeeStevenhagen}.
For the interested reader all mathematical proofs can be found in \cite{KonstKonto}.
However, in the rest of the section we will 
present a considerably simplified version of the method 
in \cite{KonstKonto}. 
%describe
%the algorithm for the construction of the polynomials with an eye to practitioner needs.
%For the interested reader all mathematical proofs can be found in \cite{KonstKonto}.
%The algorithm presented here is a considerably simplified version of the method 
%in \cite{KonstKonto}.

The Ramanujan polynomial $T_{D}(x)\in\mathbb{Z}[x]$ for $D \equiv 11 \bmod {24}$
is defined as
\begin{equation*}
  T_{D}(x) = \prod_{\tau} (x-t(\tau))
\label{tx}
\end{equation*}
for values of $\tau$ satisfying $\tau=\frac{-\beta+\sqrt{-D}}{2\alpha}$ for all primitive, reduced
quadratic forms $[\alpha, \beta, \gamma]$ of $-D$.
Every value $t(\tau)$ that corresponds to a specific form $[\alpha, \beta, \gamma]$ is defined by
\begin{equation*}
t(\tau)= (\zeta_{72}^{6k}-\zeta_{72}^{30k})\sum_{i=0}^{5} a_{2i}R_i(\tau)
\end{equation*}
where $\zeta_{72}=e^{2\pi i /72}$ and
the functions $R_i$ with $i \in \{0,1,2,3,4,5\}$ are modular functions of level $72$ and are
defined by:
$R_0(\tau)= \frac{\eta(3\tau)\eta(\tau/3)}{ \eta^2(\tau)}$,
$R_1(\tau)= \frac{\eta(3\tau)\eta(\tau/3+1/3)}{  \eta^2(\tau)}$,
$R_2(\tau)= \frac{\eta(3\tau)\eta(\tau/3+2/3)}{ \eta^2(\tau)}$,
$R_3(\tau)= \frac{\eta(\tau/3)\eta(\tau/3+2/3)}{  \eta^2(\tau)}$,
$R_4(\tau)= \frac{\eta(\tau/3)\eta(\tau/3+1/3)}{ \eta^2(\tau)}$
and $R_5(\tau)=\frac{\eta(\tau/3+2/3) \eta(\tau/3+1/3) }{ \eta^2(\tau)}$.
\REMOVED{
\begin{equation*}
R_0(\tau)= \frac{\eta(3\tau)\eta(\tau/3)}{ \eta^2(\tau)}
\end{equation*}

\begin{equation*}
R_1(\tau)= \frac{\eta(3\tau)\eta(\tau/3+1/3)}{  \eta^2(\tau)}
\end{equation*}

\begin{equation*}
R_2(\tau)= \frac{\eta(3\tau)\eta(\tau/3+2/3)}{ \eta^2(\tau)}
\end{equation*}

\begin{equation*}
R_3(\tau)= \frac{\eta(\tau/3)\eta(\tau/3+2/3)}{  \eta^2(\tau)}
\end{equation*}

\begin{equation*}
R_4(\tau)= \frac{\eta(\tau/3)\eta(\tau/3+1/3)}{ \eta^2(\tau)}
\end{equation*}

\begin{equation*}
R_5(\tau)=\frac{\eta(\tau/3+2/3) \eta(\tau/3+1/3) }{ \eta^2(\tau)}
\end{equation*}
}
The value $k$ is equal to $9\det(L_2)-8\det(L_3)$ where $\det(L_2)$ and $\det(L_3)$ are the determinants of the following
matrices $L_n$ for $n=2$ or $3$ respectively:

\begin{equation*}
L_n= \left\{
               \begin{array}{rl}
               \begin{pmatrix} \alpha & \frac{(\beta-1)}{2} \\ 0 & 1 \end{pmatrix} & \mbox{ if $n\mid\!\!\!\!/\alpha$} \\
               \begin{pmatrix} \frac{(-\beta-1)}{2} & -\gamma \\ 1 & 0 \end{pmatrix}
                                                      & \mbox{ if $n\mid \alpha$ and $n\mid\!\!\!\!/ \gamma$} \\
               \begin{pmatrix} \frac{(-\beta-1)}{2}-\alpha & \frac{(1-\beta)}{2}-\gamma \\ 1 & -1 \end{pmatrix}
                                                      & \mbox{ if $n\mid \alpha$ and $n\mid \gamma$}
               \end{array}
               \right.
\label{l2}
\end{equation*}
The values $a_{2i}$ with $i \in \{0,1,2,3,4,5\}$ are the elements
of the third row of a $6\times6$ matrix $A$. Before describing the construction of $A$ we
need to define the following two matrices:
\begin{equation*} \label{ATAS}
S_0={\begin{pmatrix}
 {0}&{\zeta_{72}^{3k}}&{0}&{0}&{0}&{0}\cr
 {0}&{0}&{\zeta_{72}^{3k}}&{0}&{0}&{0}\cr
 {\zeta_{72}^{6k}}&{0}&{0}&{0}&{0}&{0}\cr
 \\
 {0}&{0}&{0}&{0}&{{1}\over{\zeta_{72}^{3k}}}&{0}\cr
 {0}&{0}&{0}&{0}&{0}&{{1}\over{\zeta_{72}^{6k}}}\cr
 {0}&{0}&{0}&{{1}\over{\zeta_{72}^{3k}}}&{0}&{0}\cr
\end{pmatrix} },
\end{equation*}
\begin{equation*}
S_1=
{
\begin{pmatrix}
 {1}&{0}&{0}&{0}&{0}&{0}\cr
 {0}&{0}&{0}&1 \over{{\zeta_{72}^{3k}}({ {-\zeta_{72}^{30k} + \zeta_{72}^{6k}} }})&{0}&{0}\cr
 {0}&{0}&{0}&{0}&{{{\zeta_{72}^{3k}}\over{-\zeta_{72}^{30k} + \zeta_{72}^{6k}}}}&{0}\cr
 {0}&{ {\zeta_{72}^{3k}}({-\zeta_{72}^{30k} + \zeta_{72}^{6k}}) }&{0}&{0}&{0}&{0}\cr
 {0}&{0}&{{-\zeta_{72}^{30k} + \zeta_{72}^{6k}}\over{\zeta_{72}^{3k}}}&{0}&{0}&{0}\cr
 {0}&{0}&{0}&{0}&{0}&{1}\cr
 \end{pmatrix} }.
\end{equation*}
Using $S_0$ and $S_1$ we can compute four new matrices $T_2 = S_0^{9}$, $T_3 = S_0^{-8}$, $S_2 = S_0^{-1} S_1 S_0^{-10} S_1 S_0^{-1} S_1 S_0^{-18}$ and
$S_3 = S_0^{-1} S_1 S_0^{7} S_1 S_0^{-1} S_1 S_0^{16}$.
\REMOVED{
\begin{equation}
S_2 = S_0^{-1} S_1 S_0^{-10} S_1 S_0^{-1} S_1 S_0^{-18}
\end{equation}

\begin{equation}
S_3 = S_0^{-1} S_1 S_0^{7} S_1 S_0^{-1} S_1 S_0^{16}
\end{equation}

\begin{equation}
T_2 = S_0^{9}
\end{equation}
and
\begin{equation}
T_3 = S_0^{-8}
\end{equation}
}
Now the matrix $A$ is equal to $A_{2} A_3 B$ where $B$ is equal to
\begin{equation*}
B=\left\{
               \begin{array}{rl}
               \begin{pmatrix}
1 & 0 & 0 & 0& 0 & 0 \\
0 & \zeta_{72}^{k-1} & 0 & 0 & 0& 0 \\
0 & 0 & \zeta_{72}^{2k -2} & 0 & 0 & 0 \\
0 & 0 & 0 & \zeta_{72}^{2k-2} & 0 & 0 \\
0 & 0 & 0 & 0 & \zeta_{72}^{k-1} &0 \\
0 & 0 & 0 & 0 & 0 & \zeta_{72}^{3k-3}
\end{pmatrix}  \mbox{ if } k\equiv 1 \bmod 3 \\ \\ \\

\begin{pmatrix}
1 & 0 & 0 & 0& 0 & 0 \\
0 & 0 &  \zeta_{72}^{k-2}  & 0 & 0& 0 \\
0 & \zeta_{72}^{2k -1} &0 & 0 & 0 & 0 \\
0 & 0 & 0 & 0 & \zeta_{72}^{2k-1}  & 0 \\
0 & 0 & 0  & \zeta_{72}^{k-2} &0  &0 \\
0 & 0 & 0 & 0 & 0 & \zeta_{72}^{3k-3}
\end{pmatrix}  \mbox{ if } k\equiv 2 \bmod 3

                            \end{array}
                            \right.
\end{equation*}
and
\begin{equation*}
A_n= \left\{
  \begin{array}{ll}
                S_nT_n^{\frac{1}{\alpha} \bmod {N(n)}}S_nT_n^{-\alpha}S_nT_n^{(\frac{1}{\alpha}(\frac{\beta-1}{2}) -1)  \bmod {N(n)}} & \mbox{ if $n\nmid \alpha$} \\
                T_n^{(1-\frac{\beta+1}{2}) \bmod {N(n)}}S_nT_nS_nT_n^\gamma  & \mbox{ if $n\mid \alpha$ and $n\nmid \gamma$} \\
                T_n^{(1-\frac{\beta+1}{2}-\alpha) \bmod {N(n)}}S_nT_nS_nT_n^{(-1+\alpha+\beta+\gamma) \bmod N(n)} & \mbox{ if $n\mid \alpha$ and $n\mid \gamma$}
               \end{array}
               \right.
\label{Gmatrix}
\end{equation*}
for $n=2,3$ and $N(2)=8$,$N(3)=9$.
\REMOVED{
\begin{equation}
F= \left\{
               \begin{array}{rl}
                S_1T_1^{(\frac{1}{\alpha} \bmod 8)}S_1T_1^{-\alpha}S_1T_1^{(\frac{1}{\alpha}((\frac{1}{\alpha} \bmod 8)(\frac{\beta-1}{2}) -1) )} & \mbox{ if $2\mid\!\!\!\!/\alpha$} \\
                T_1^{(1-\frac{\beta+1}{2}) \bmod 8}S_1T_1S_1T_1  & \mbox{ if $2\mid \alpha$ and $2\mid\!\!\!\!/ \gamma$} \\
                T_1^{(1-\frac{\beta+1}{2}-\alpha) \bmod 8}S_1T_1S_1T_1^{(\alpha+\beta+\gamma)^{-1}((\alpha+\frac{\beta+1}{2}-1) \bmod 8+\gamma + \frac{\beta-1}{2})} & \mbox{ if $2\mid \alpha$ and $2\mid \gamma$}
               \end{array}
               \right.
\label{Fmatrix}
\end{equation}
}
\REMOVED{
and
\begin{equation*}
B=\left\{
               \begin{array}{rl}
               \begin{pmatrix}
1 & 0 & 0 & 0& 0 & 0 \\
0 & \zeta_{72}^{k-1} & 0 & 0 & 0& 0 \\
0 & 0 & \zeta_{72}^{2k -2} & 0 & 0 & 0 \\
0 & 0 & 0 & \zeta_{72}^{2k-2} & 0 & 0 \\
0 & 0 & 0 & 0 & \zeta_{72}^{k-1} &0 \\
0 & 0 & 0 & 0 & 0 & \zeta_{72}^{3k-3}
\end{pmatrix}  \mbox{ if } k\equiv 1 \bmod 3 \\ \\ \\

\begin{pmatrix}
1 & 0 & 0 & 0& 0 & 0 \\
0 & 0 &  \zeta_{72}^{k-2}  & 0 & 0& 0 \\
0 & \zeta_{72}^{2k -1} &0 & 0 & 0 & 0 \\
0 & 0 & 0 & 0 & \zeta_{72}^{2k-1}  & 0 \\
0 & 0 & 0  & \zeta_{72}^{k-2} &0  &0 \\
0 & 0 & 0 & 0 & 0 & \zeta_{72}^{3k-3}
\end{pmatrix}  \mbox{ if } k\equiv 2 \bmod 3

                            \end{array}
                            \right.
\end{equation*}
}

It is easy to see that  every row in the matrix $A$ has only one non zero element. Thus, only one value $a_{2i}$ is not equal to zero and the computation of every value $t(\tau)$ requires the evaluation of only one value $R_i(\tau)$.

\subsection{Transformation of the Roots}

In order to use Ramanujan polyomials in the CM method, we must prove that they have roots modulo $p$ and then
find a transformation of their
roots modulo $p$ to the roots modulo $p$ of the corresponding Hilbert polynomials. The following theorem
proves that a Ramanujan polynomial with degree $h$ has exactly $h$ roots modulo $p$
under certain conditions (which are satisfied in the CM method):
\begin{theorem}
A Ramanujan polynomial $T_D(x)$ with degree $h$ has exactly $h$ roots modulo $p$ if and only if
the equation $4p = u^2 + Dv^2$ has integer solutions and $p$ does not divide the
discriminant $\Delta(T_{D})$ of the polynomial.
\label{ramanujanroots}
\end{theorem}
\begin{proof}
Let
$H_K$ be the Hilbert class field of the imaginary  quadratic field
$K=\mathbb{Q}(\sqrt{-D})$, and
let $\mathcal{O}_{H_K}$ and $\mathcal{O}_K$ be the rings of algebraic integers
of $H_K$ and $K$ respectively.
Let $p$ be a prime such that $4p=u^2+Dv^2$ has integer solutions. Then,
according to \cite[Th. 5.26]{C89} $p$ splits completely in $H_K$.
Proposition 5.29 in
\cite{C89} implies that (since $t_D$ generates $H_K$) $T_D(x)$ has a root modulo $p$ if and only
if $p$ splits in $H_K$ and does not divide its
discriminant
%\footnote{For a definition of the discriminant of a
%polynomial see \cite{C93}.} 
$\Delta(T_{D})$. But since
$\frac{\mathcal{O}_{H_K}}{p \mathcal{O}_{H_K}}/ \mathbb{F}_p$ is
Galois, $T_D(x)$ has not only one root modulo $p$, but $h$
distinct roots modulo $p$.
\end{proof}

We will present now a method to
retrieve a root modulo $p$ of the Hilbert polynomial $H_{D}(x)$
from a root
modulo $p$ of the corresponding Ramanujan polynomial $T_{D}(x)$.
Our aim is to find a transformation that maps a real root of the Ramanujan polynomial to a real root of the
corresponding Hilbert polynomial. Then, we can reduce this transformation
modulo a prime ideal  of the ring of integers of the Hilbert class field.
In this way we see that the same transformation will transfer a root of the Ramanujan polynomial modulo $p$
to a root of the Hilbert polynomial modulo $p$.
%and it will give us a method to compute the $j$-invariant
%and construct the desired elliptic curve.
We know that if $\ell_0=(1,1,\frac{1+D}{4})$ is a quadratic form (known as the principal form) that corresponds
to the root $\tau_{\ell_0}=-\frac{1}{2}+i \frac{\sqrt{-D}}{2}$ then $j(\tau_{\ell_0})$ is a real root of the
Hilbert polynomial $H_D(x)$. The following lemma shows that $t_D$ is a real root of the Ramanujan polynomial
$T_D(x)$.

\begin{lemma} \label{modtn}
The value $t_D $ is a real root of the Ramanujan polynomial $T_D(x)$ and is equal to:
\[
t_D=\sqrt{3} R_2(\tau_{\ell_0}).
\]
\end{lemma}

\begin{proof}
Set
\begin{equation*}
q_D=\exp(-\pi \sqrt{D})=-\exp(2\pi i \tau_{\ell_0}),
\end{equation*}
where $\tau_{\ell_0}=-\frac{1}{2} + i \frac{\sqrt{-D}}{2}$.
Then
\begin{equation*}
f(q_D)= f( -\exp(2\pi i \tau_{\ell_0}))=  \exp(2\pi i \tau_{\ell_0})^{-1/24} \eta(\tau_{\ell_0}),
\end{equation*}
\begin{equation*}
f(q_D^3)= \exp(2\pi i \tau_{\ell_0})^{-3/24} \eta(3\tau_{\ell_0}),
\end{equation*}
\begin{equation*}
f(q_D^{1/3})=\exp(2\pi i \tau_{\ell_0})^{-\frac{1}{3\cdot 24}}\eta(\frac{\tau_{\ell_0}}{3}).
\end{equation*}
Taking Eq.~(\ref{tndef}) and all the above equations into consideration
we can easily derive that $t_D=\sqrt{3} R_2(\tau_{\ell_0})$.

If we could prove that $t(\tau_{\ell_0}) = \sqrt{3} R_2(\tau_{\ell_0})$ then it will immediately follow
that $t_D = t(\tau_{\ell_0})$ and thus it is a root of the Ramanujan polynomial.
We have that
\[
 t(\tau_{\ell_0})=(\zeta_{72}^6 -\zeta_{72}^{30}) R_2(\tau_{\ell_0}),
\]
since $k=1$ and the matrix $A=A_2 A_3 B$ is by computation equal to the identity matrix for every discriminant $D$. Notice  that
the principal form equals $[\alpha,\beta,\gamma]=[1,1,\frac{1-D}{4}]$, therefore $2,3 \nmid \alpha=1$ and
$L_2=L_3=\mathrm{Id}_2$, $B=\mathrm{Id}_6$ and $A_n=S_nT_n^{\frac{1}{\alpha} \bmod {N(n)}}S_nT_n^{-\alpha}S_nT_n^{(\frac{1}{\alpha}(\frac{\beta-1}{2}) -1)  \bmod {N(n)}}$ for $n=2,3$.
Finally, observe that $\sqrt{3}=\zeta_{72}^{6} -\zeta_{72}^{30}.$  Indeed, the value
 $i \sqrt{3}$ can be expressed  as a difference of two primitive $3$-roots of unity $\zeta_3,\zeta_3^2$ since
 $i=\zeta_{72}^{18}$  and $\zeta_3=\zeta_{72}^{24}$. Thus $t(\tau_{\ell_0}) = \sqrt{3} R_2(\tau_{\ell_0}) = t_D$.
\end{proof}

%We will find a transformation $T(t_D)=j(\tau_{\ell_0})$  of a real root $t_D$ of the Ramanujan polynomial
%to a real root $j(\tau_{\ell_0})$ of the corresponding Hilbert polynomial.

%
% Thus, the values $t_D$ correspond to real roots of the Ramanujan polynomial
% (it is known that when $\tau^*$ corresponds to a principal form, then any class
% invariant $g(\tau^*)$ of a class polynomial is a real root {\bf auto einai swsto?}).

%It is also known that when $\tau^*$ corresponds to a principal form, then the class
%invariant $j(\tau^*)$ of $H_{D}(x)$ is a real root.
%Hence the goal is to find a transformation
%$T(~)$ such that $j(\tau_{\ell_0})=T(t_D)$.
%The following lemma presents the necessary transformation.
%
\begin{lemma}
\label{lemma_trans}
Suppose $R_T$ is a real root of a Ramanujan polynomial $T_D(x)$.
Then, the real number $R_H$ obtained from the equation 
\begin{equation}
\label{eq_trans}
R_H=(R_T^6-27R_T^{-6}-6)^3
\end{equation}
is a real root of the corresponding Hilbert polynomial $H_{D}(x)$. 
\end{lemma}

\begin{proof}
Set $R_T = t_D$ and $R_H = j(\tau_{\ell_0})$.
Using Equations (4.4) and (4.5) from \cite{Berndt-Chan} it can be easily derived that
%\begin{equation*}
$h(e^{2\pi i\tau_{\ell_0}/3})-27h(e^{2\pi i\tau_{\ell_0}/3})^{-1}=\gamma_2(\tau_{\ell_0})+6$
%\end{equation*}
where $\gamma_2^3(\tau_{\ell_0})=j(\tau_{\ell_0})$ and
\begin{equation}
\label{h_eq}
h(q)= \frac{f^{12}(-q^3)}{qf^6(-q)f^6(-q^9)}.
\end{equation}
%
%and
%\begin{equation*}
%$\gamma_2^3(\tau_{\ell_0})=j(\tau_{\ell_0})$.
%\end{equation*}
Thus, $j(\tau_{\ell_0})=(h(e^{2\pi i\tau_{\ell_0}/3})-27h(e^{2\pi i\tau_{\ell_0}/3})^{-1}-6)^3$ which means that
we now have to find the relation between $t_D$ and $h(e^{2\pi i\tau_{\ell_0}/3})$.
Substituting $q$ with $e^{2\pi i\tau_{\ell_0}/3}$ in Eq.~(\ref{h_eq}) we have that
$h(e^{2\pi i\tau_{\ell_0}/3})=\frac{f^{12}(-e^{2\pi i \tau_{\ell_0}})}{e^{2\pi i\tau_{\ell_0}/3}f^6(-e^{2\pi i\tau_{\ell_0}/3})f^6(-e^{3(2\pi i\tau_{\ell_0})})}$. 
Noticing that $q_D=\exp(-\pi \sqrt{D})=-\exp(2\pi i \tau_{\ell_0})$ and from Eq.~(\ref{tndef}) %and (\ref{qndef})
we derive that
$h(e^{2\pi i\tau_{\ell_0}/3})=-27t_D^{-6}$ which completes the proof of the lemma.
\end{proof}

The final step is to reduce Eq.~(\ref{eq_trans}) modulo $p$. The elements  $R_H,R_T$
are not in $\mathbb{Z}$ but are elements of the ring of algebraic integers  $\mathcal{O}_{H_K}$ of the Hilbert class field and
can be reduced modulo an ideal $P$ extending the ideal $p\mathbb{Z}$ of $\mathbb{Z}$.
But the ideal $p\mathbb{Z}$ splits completely, therefore the Galois extension $\frac{\mathcal{O}_{H_K}/P}{\mathbb{Z}/p\mathbb{Z} }$ is
the trivial one, and $\mathcal{O}_{H_K}/P$ is the field $\mathbb{F}_p$.
The argument above proves that Eq.~(\ref{eq_trans}) holds not only for the real roots of the polynomials
but also for their roots modulo $p$.
The interested reader is referred to \cite{C89,S04,ST87} for definitions on
algebraic number theory
not given here. 
%due to lack of space.
%
%
Using Eq.~(\ref{eq_trans}), we can easily derive the modular polynomial $\Phi_T(x,j)$ for Ramanujan polynomials. 
The polynomial will be equal to:
\begin{equation}
\label{ramanujan-phi-eq}
\Phi_T(x,j) = (x^{12}-6x^6-27)^3 -jx^{18}.
\end{equation}
%Thus, following the estimation in \cite{EM02}, we conclude 
%that the precision required for the construction of the Ramanujan polynomials is approximately
%$\frac{1}{36}\frac{\pi\sqrt{D}}{\ln 2} \sum_{\tau} \frac{1}{\alpha}$. 

\section{Precision Requirements for the Construction of the Polynomials}
\label{prec-section}

In this section we focus on the precision required for the construction
of all previously mentioned polynomials. 
In order to compare them, we introduce the notion of {\em logarithmic height} for 
estimating the size of a polynomial. For a polynomial $g(x)=\sum_{i=0}^n a_i x^i \in \mathbb{Z}[x]$ its 
logarithmic height is  defined as
\[
 H(g)=\max_{i=0,\ldots,n} \log_2 |a_i|.
\]
The value $H(g)$ is actually the bit-precision needed for performing all floating point computations 
in order to obtain the coefficients
of the polynomial $g(x)$. 

Starting from Hilbert polynomials, an estimation of their precision requirements in bits (and of their
logarithmic height also) was
given in \cite{LZ94}:
\begin{equation*}
\mbox{H-Prec}(D) \approx \frac{\ln 10}{\ln 2} (h/4 + 5) +
\frac{\pi\sqrt{D}}{\ln 2} \sum_{\tau} \frac{1}{\alpha}
\label{hplz}
\end{equation*}
with the sum running over the same values of $\tau$ as the product
in Eq.~(\ref{hx}). 
A slightly different bound was given in \cite{M90} which is remarkably accurate:
\begin{equation*}
\mbox{H-Prec1}(D) \approx 33 +
\frac{\pi\sqrt{D}}{\ln 2} \sum_{\tau} \frac{1}{\alpha}.
\label{hplz1}
\end{equation*}
It will be shown in the rest of the section that based on this estimation, we can derive estimations of the
precision requirements of every class polynomial.

Let $f$ be a modular function, such that $f(\tau)$ for some $\tau \in \mathbb{Q}(\sqrt{-D})$ generates the Hilbert class field of 
$\mathbb{Q}(\sqrt{-D})$.
The element $f(\tau)$ is an algebraic integer, and let us denote by $P_f$ its minimal polynomial. 
For every modular function there is a polynomial $\Phi$ (called modular polynomial) such that
$\Phi(f,j) = 0$ where $j$ is the modular function used in the construction of Hilbert polynomials.
This polynomial equation is used (as we show in the previous section) in order to transform the roots of 
the minimal polynomial 
of a class invariant to the roots of the Hilbert polynomial. We have seen that in the cases 
of Weber, $M_{D,l}(x)$ and Ramanujan polynomials the degree in $j$ of the modular polynomial is equal to 1
while for $M_{D,p_1, p_2}(x)$ polynomials is at least 2.
Asymptotically, one can estimate the ratio of the logarithmic height $h(j(\tau))$ of the algebraic integer $j(\tau)$
to the logarithmic height $h(f(\tau))$ of the algebraic integer $f(\tau)$
\footnote{Let $K$ be a number field, $\alpha\in K$ be an algebraic
number and $M_K$ be the set of absolute values on $K$. Following
the notation of \cite[VIII]{S86}, the absolute logarithmic height
of an element $\alpha \in K$ is defined as $h(\alpha) =
\frac{1}{[K:\mathbb{Q}]} \log_2 \left( \prod_{v\in M_K} \max \{
|\alpha|_v, 1 \} \right)$. }
. Namely, 
\begin{equation}
\label{ratio11}
 \lim_{h(j(\tau))\rightarrow \infty} \frac{h(j(\tau))}{h(f(\tau))}=\frac{\deg_f \Phi(f,j)}{\deg_j \Phi(f,j)}=r(f),
\end{equation}
where the limit is taken over all CM-points $\mathrm{SL}_2(\mathbb{Z})\tau \in \mathbb{H}$ \cite{HinSil}.
Concerning Weber polynomials, we can easily compute the values of $r(f)$ from 
Eq. ~(\ref{eq-phi-1}) and Eq.~(\ref{eq-phi-2}).
Thus, when $D \not\equiv 0 \pmod 3$, $r(f)$ will be equal to 24 and when $D \equiv 0 \pmod 3$, $r(f)$ will be equal to 8.

A question that immediately arises is how Eq.~(\ref{ratio11}) can be used for the estimation of the logarithmic
height of the minimal polynomial $P_f$. The following Lemma gives an answer to this question.
\begin{lemma}
Suppose that $H(P_f)$ is the logarithmic height of the minimal polynomial of the algebraic number $f(\tau)$ and
$H(P_j)$ is the logarithmic height of the corresponding Hilbert polynomial. If $f(\tau)$ generates the Hilbert 
class field then
\begin{equation}
\label{ratio12}
 \lim_{h(j(\tau))\rightarrow \infty} \frac{H(P_j)}{H(P_f)}=\frac{\deg_f \Phi(f,j)}{\deg_j \Phi(f,j)}=r(f).
\end{equation}
If $f(\tau)$ does not generate the Hilbert 
class field but an algebraic extension of it with extension degree $m$ then
\begin{equation*}
%\label{ratio13}
 \lim_{h(j(\tau))\rightarrow \infty} \frac{H(P_j)}{H(P_f)}=\frac{\deg_f \Phi(f,j)}{\deg_j \Phi(f,j)}=\frac{r(f)}{m}.
\end{equation*}

\end{lemma}
\begin{proof}
The proof is based on the following bounds\cite[Th. 5.9]{S86}:
\begin{equation*} 
%\label{roots-coeff}
 -k + k h(a) \leq H(P_a) \leq k-1+ k h(a)
\end{equation*}
where $h(a)$ is the logarithmic height of the algebraic integer $a$ and $k$ is the degree of its minimal polynomial $P_a$. 
If $f(\tau)$ generates the Hilbert class field then the degree of its minimal polynomial is equal to the degree of the
corresponding Hilbert polynomial. Suppose that their degree is equal to $k$. Then, we have that
\begin{equation}
\label{eq-f} 
 -k + k h(f(\tau)) \leq H(P_f) \leq k-1+ k h(f(\tau))
\end{equation}
and
\begin{equation*} 
 -k + k h(j(\tau)) \leq H(P_j) \leq k-1+ k h(j(\tau)).
\end{equation*}
Thus,
\begin{equation*} 
\frac{-k + k h(j(\tau))}{k-1+ k h(f(\tau))} \leq \frac{H(P_j)}{H(P_f)} \leq \frac{k-1+ k h(j(\tau))}{-k + k h(f(\tau))}.
\end{equation*}
Taking the limit $h(j(\tau))\rightarrow \infty$ we obtain that
\begin{equation}
\label{eq-degree}
\frac{H(P_j)}{H(P_f)}\rightarrow r(f).
\end{equation}
In the case that $f(\tau)$ generates an algebraic extension of the Hilbert class field, we similarly have that
\begin{equation}
\label{eq-ext-degree}
\frac{H(P_j)}{H(P_f)}\rightarrow \frac{r(f)}{m}
\end{equation}
where $m$ is the degree of the extension. This is easily derived from the fact that the degree of the minimal polynomial
$P_f$ is $m$ times larger than the degree of the corresponding Hilbert polynomial and Eq.~(\ref{eq-f}) becomes
\begin{equation*}
%\label{eq-f1} 
 -mk + mk h(f(\tau)) \leq H(P_f) \leq mk-1+ mk h(f(\tau)).
\end{equation*}
Thus,
\begin{equation*} 
\frac{-k + k h(j(\tau))}{mk-1+ mk h(f(\tau))} \leq \frac{H(P_j)}{H(P_f)} \leq \frac{k-1+ k h(j(\tau))}{-mk + mk h(f(\tau))}.
\end{equation*}
\end{proof}

Eq.~(\ref{eq-degree}) and Eq.~(\ref{eq-ext-degree}) 
relate the precision required for the construction of Hilbert polynomials with the precision
needed for other classes of polynomials. 
Estimating the height $H(P_j)$ of Hilbert polynomials with the quantity $\frac{\pi\sqrt{D}}{\ln 2} \sum_{\tau} \frac{1}{\alpha}$, 
we can derive the precision requirements for the construction of every class polynomial by the equation:
\begin{equation*}
\frac{m}{r(f)} \frac{\pi\sqrt{D}}{\ln 2} \sum_{\tau} \frac{1}{\alpha},
%\label{prec-est-morain}
\end{equation*}
where $m$ is either 1 or larger.

Obviously,  
we want to find class invariants $f(\tau)$ so that the ratio  $r(f)$ is as big as possible.
However, there is a limit on the ratio $r(f)$.
It is known \cite{BS06} that $r(f) \leq 800/7$ and if the Selberg eigenvalue conjecture in \cite{Sarnak} holds
then $r(f) \leq 96$. 
Concerning Weber polynomials, 
%we can easily compute the values of $r(f)$ from 
%Tables~\ref{trans:not-mod-3} and ~\ref{trans:mod-3}.
%For example, when $D \equiv 7 \pmod 8 \not\equiv 0 \pmod 3$ the precision requirements of the 
%corresponding Weber polynomials with be approximately $\frac{1}{72}\frac{\pi\sqrt{D}}{\ln 2} \sum_{\tau} \frac{1}{\alpha}$.
%Similarly, we can estimate the precision requirements of Weber polynomials for all cases of $D$. However, 
when $D \equiv 3 \pmod 8$ their degree
is three times larger than the degree of the corresponding Hilbert polynomials. Therefore, for this case of
$D$, the estimation of the precision requirements will be approximately
$\frac{3}{r(f)}\frac{\pi\sqrt{D}}{\ln 2} \sum_{\tau} \frac{1}{\alpha}$. 
Concluding, an estimation of the precision requirements of Weber polynomials will be equal to 
$\frac{1}{24}\frac{\pi\sqrt{D}}{\ln 2} \sum_{\tau} \frac{1}{\alpha}$ for $D \not\equiv 0 \pmod 3$
and $\frac{1}{8}\frac{\pi\sqrt{D}}{\ln 2} \sum_{\tau} \frac{1}{\alpha}$ for $D \equiv 0 \pmod 3$.

Based again on Eq.~(\ref{ratio12}), it can be concluded that 
the precision required for the construction of the
$M_{D, l}(x)$ polynomials is approximately $\frac{1}{(l+1)}\frac{\pi\sqrt{D}}{\ln 2} \sum_{\tau} \frac{1}{\alpha}$
and for
$M_{D, p_1, p_2}(x)$ polynomials is approximately $\frac{(p_1-1)(p_2-1)}{12(p_1+1)(p_2+1)}\frac{\pi\sqrt{D}}{\ln 2} \sum_{\tau} \frac{1}{\alpha}$
where the sum runs over the same values of $\tau$ as the product
in Eq.~(\ref{hx}) \cite{EM02}.
Thus, it is equal to $\frac{1}{28}\frac{\pi\sqrt{D}}{\ln 2} \sum_{\tau} \frac{1}{\alpha}$ for $M_{D, 3, 13}(x)$ polynomials and to $\frac{1}{24}\frac{\pi\sqrt{D}}{\ln 2} \sum_{\tau} \frac{1}{\alpha}$ for $M_{D, 5, 7}(x)$ polynomials. The above precision estimations are summarized in Table~\ref{prec-ml}.
\begin{table}
\centering
\begin{tabular}{||l|l||} \hline
 class polynomial & precision estimation \\ \hline \hline
 $M_{D,3}(x)$ & $\frac{1}{4}\frac{\pi\sqrt{D}}{\ln 2} \sum_{\tau} \frac{1}{\alpha}$ \\ \hline
 $M_{D,5}(x)$ & $\frac{1}{6}\frac{\pi\sqrt{D}}{\ln 2} \sum_{\tau} \frac{1}{\alpha}$ \\ \hline
 $M_{D,7}(x)$ & $\frac{1}{8}\frac{\pi\sqrt{D}}{\ln 2} \sum_{\tau} \frac{1}{\alpha}$ \\ \hline
 $M_{D,13}(x)$ & $\frac{1}{14}\frac{\pi\sqrt{D}}{\ln 2} \sum_{\tau} \frac{1}{\alpha}$ \\ \hline
 $M_{D,5,7}(x)$ & $\frac{1}{24}\frac{\pi\sqrt{D}}{\ln 2} \sum_{\tau} \frac{1}{\alpha}$ \\ \hline
 $M_{D,3,13}(x)$ & $\frac{1}{28}\frac{\pi\sqrt{D}}{\ln 2} \sum_{\tau} \frac{1}{\alpha}$ \\ \hline
\end{tabular}
\caption{Precision estimations for $M_{D,l}(x)$ and $M_{D,p_1,p_2}(x)$ polynomials.}
\label{prec-ml}
\end{table}
Finally, in order to find an estimation for the precision requirements of Ramanujan polynomials, we use
 Eq.~(\ref{ratio12}) and Eq.~(\ref{ramanujan-phi-eq}). 
%But first we must derive the modular polynomial $\Phi_T(x,j)$ for Ramanujan polynomials.
%Using Eq.~(\ref{eq_trans}), we find that the modular polynomial is equal to:
%\begin{equation*}
%\Phi_T(x,j) = (x^{12}-6x^6-27)^3 -jx^{18}.
%\end{equation*}
%Thus, following Eq.~(\ref{ratio12}), 
We easily conclude 
that the precision required for the construction of the Ramanujan polynomials is approximately
$\frac{1}{36}\frac{\pi\sqrt{D}}{\ln 2} \sum_{\tau} \frac{1}{\alpha}$.

%In conclusion, the ranking of the polynomials based on the precision requirements (starting from the smallest) is
%$W_{D \equiv 7 \pmod 8 \not\equiv 0 \pmod 3}(x)$ $<$ $T_D(x)$ = $W_{D/4 \equiv 2,6 \pmod 8 \not\equiv 0 \pmod 3}(x)$ =
%$W_{D/4 \equiv 1 \pmod 8 \not\equiv 0 \pmod 3}(x)$ $<$ $M_{D,3,13}(x)$ $<$ $M_{D,5,7}(x)$ =
%$W_{D \equiv 7 \pmod 8 \equiv 0 \pmod 3}(x)$ = $W_{D \equiv 3 \pmod 8 \not\equiv 0 \pmod 3}(x)$ $<$
%$W_{D/4 \equiv 5 \pmod 8 \not\equiv 0 \pmod 3}(x)$ $<$ $M_{D,13}(x)$ $<$ $W_{D/4 \equiv 2,6 \pmod 8 \equiv 0 \pmod 3}(x)$ =
%$W_{D/4 \equiv 1 \pmod 8 \equiv 0 \pmod 3}(x)$ $<$ $W_{D \equiv 3 \pmod 8 \equiv 0 \pmod 3}(x)$ =
%$M_{D,7}(x)$ $<$ $M_{D,5}(x)$ = $W_{D/4 \equiv 5 \pmod 8 \equiv 0 \pmod 3}(x)$ $<$ $M_{D,3}(x)$.   
%However, the estimations of the precision requirements presented above, are asymptotic. When moderate values of discriminant
%$D$ are used, these estimations slightly change and so does the ranking of the polynomials.

\section{Implementation and Experimental Results}
\label{exper}

In this section, we discuss some issues regarding the
construction of the Weber, $M_{D,l}(x)$, $M_{D,p_1,p_2}(x)$ and Ramanujan polynomials. All
implementations and experiments were made in Pari 2.3.1  \cite{PARI2} compiled with GMP-4.2.1 kernel  \cite{gnu:gnu}  and have been carried out on a 
double 2GHz Xeon machine   running Linux 2.6.9-22  and equipped with 2Gb of main memory.

%implementation of the Complex Multiplication method for the generation of prime order ECs
%and our experimental results concerning its time and space efficiency.
%All of our implementations were made in ANSI C using the (ANSI C)
%GNUMP~\cite{gnu:gnu} library for high precision floating point
%arithmetic and also for the generation and manipulation of
%integers of unlimited precision. The implementation focuses mainly on the
%construction of the Weber, $M_{D,l}(x)$, $M_{D,p_1,p_2}(x)$ and Ramanujan polynomials. All
%implementations and experiments have been carried out on a Pentium
%4 (2.80 GHz) running Linux and equipped with 504 MB of main
%memory.
 
 \begin{figure}
 \begin{minipage}{7cm}
 \centering\epsfig{file=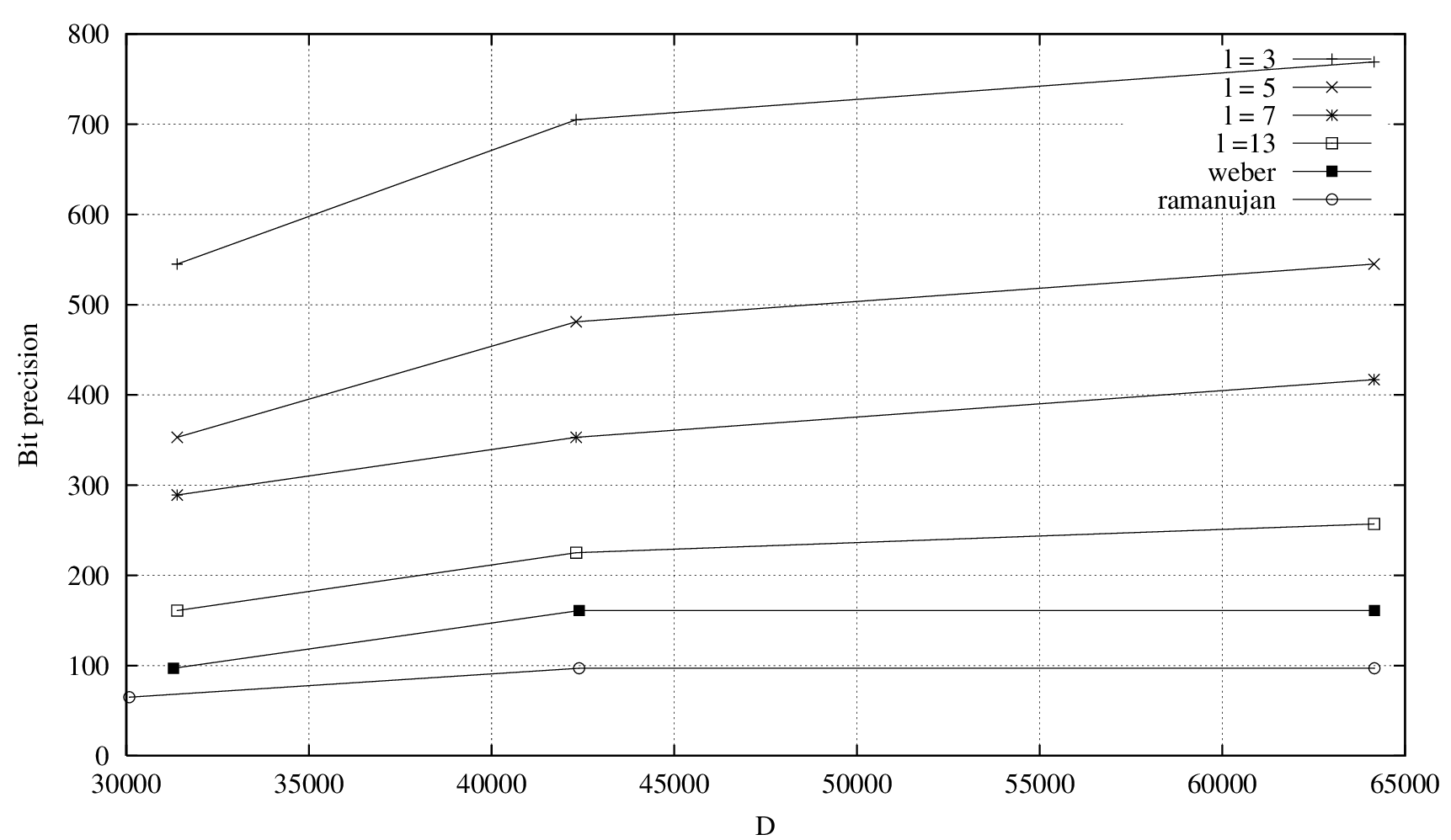, width=4.7cm, angle = 270}
 \end{minipage} \hfill
 \begin{minipage}{7cm}
 \centering\epsfig{file=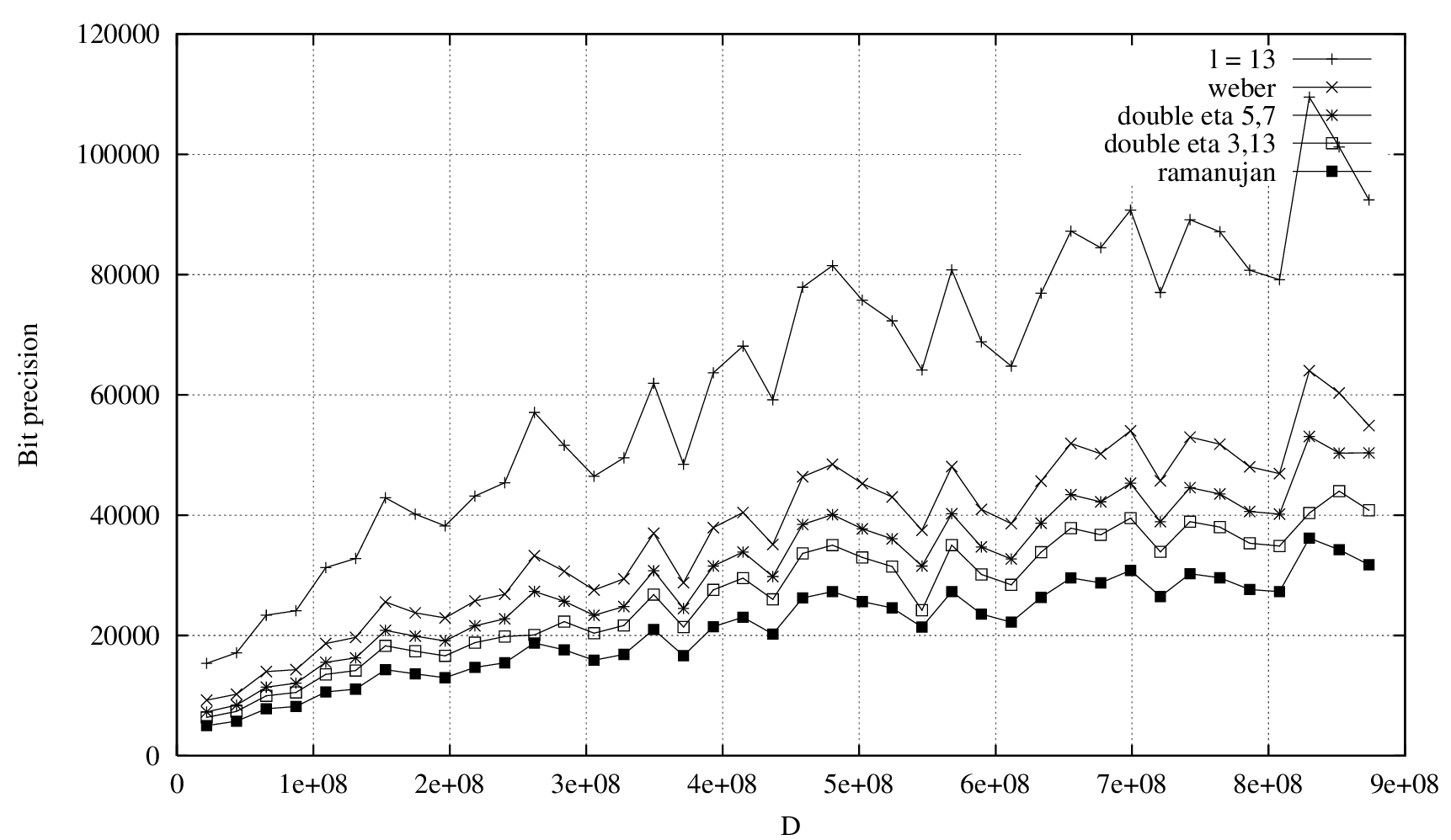, width=4.7cm, angle = 270}
 \end{minipage}
 \caption{Bit precision for the construction of
 all polynomials.}
 \label{hwprecision}
 \end{figure}

In Figure~\ref{hwprecision} we report on the precision needed for the
construction of all polynomials for various values of $D$.
In the left figure, we examine the precision requirements of Ramanujan, Weber ($D \not\equiv 0 \pmod 3$) and
$M_{D,l}(x)$ polynomials for all values of $l$.
The values of $D$ range from 30083 to 64163 while the degree $h$ ranges from 32 to 48.
We noticed, as the
theory dictates, that the precision required for the construction of Ramanujan polynomials is much less than the precision required for the construction of Weber and
$M_{D,l}(x)$ polynomials for all
values of $D$ that we examined. Weber polynomials require less precision than $M_{D,l}(x)$ polynomials, while
among them $M_{D,13}(x)$ polynomials require the least precision. 
Examining larger values of the discriminant $D$ and adding $M_{D,3,13}(x)$ and
$M_{D,5,7}(x)$ polynomials in our comparison, we show (Figure~\ref{hwprecision} (right)) that
Ramanujan polynomials are constructed more efficiently than all other polynomials. $M_{D,3,13}(x)$ polynomials require less precision
than $M_{D,5,7}(x)$ polynomials which are constructed more efficiently than Weber polynomials.
In this figure, we examined all values of $D$ from 21840299 to 873600299 using a step of 21840000. 
The degree $h$ of
the constructed polynomials (for these values of $D$) ranges from 2880 to 17472.
Summarising the results of our experiments, we see that Ramanujan polynomials outweight
$M_{D,13}(x)$, Weber, $M_{D,5,7}(x)$ and $M_{D,3,13}(x)$ polynomials as they require on average 
66\%, 42\%, 32\%
and 22\% less precision respectively.
Table~\ref{prectable} shows this difference by presenting the exact bit precision needed for
the construction of the polynomials for several values of $D$.
%45\%, 69\%, 34\%
%and 25\% less precision respectively ??? auta einai swsta??? bres to m.o. tou ratio.

\begin{table}[h]
\begin{center}
\begin{tabular}{|p{1.8cm}|p{1.0cm}|p{1.8cm}|p{1.8cm}|p{1.8cm}|p{1.8cm}|p{1.8cm}|} \hline
 $D$   & $h$  & $M_{D,13}(x)$ & Weber & $M_{D,5,7}(x)$ & $M_{D,3,13}(x)$ & Ramanujan \\ \hline \hline
 109200299  & 5016  & 31270  & 18657  & 15546 & 13534  & 10624   \\ \hline
 240240299  & 6944  & 45402  & 26837  & 22757 & 19834  & 15442   \\ \hline
 349440299  & 9772  &  61933 & 37004  & 30768 &  26804 &  20998  \\ \hline
 458640299  & 12660  &  77894 &  46387 & 38447 &  33633 & 26245   \\ \hline
 698880299  & 13950  &  90734 &  54030 & 45311 &  39508 & 30813  \\ \hline
 851760299  & 15904  &  101214 &  60333 & 50322 & 43984 & 34243  \\ \hline
\end{tabular}
\end{center}
\caption{Precision requirements (in bits) for the computation of $M_{D,13}(x)$, Weber, $M_{D,5,7}(x)$, $M_{D,3, 13}(x)$ and Ramanujan polynomials.}
\label{prectable}
\end{table}

\begin{table}[h]
\begin{center}
\begin{tabular}{|p{1.8cm}|p{1.0cm}|p{1.8cm}|p{1.8cm}|p{1.8cm}|p{1.8cm}|p{1.8cm}|} \hline
 $D$   & $h$  & $M_{D,13}(x)$ & Weber & $M_{D,5,7}(x)$ & $M_{D,3,13}(x)$ & Ramanujan \\ \hline \hline
 109200299  & 5016  &  134 &  245  &   68 &    59  &  47     \\ \hline
 240240299  & 6944  &  271  & 492 &  138 &  119  &  94     \\ \hline
 349440299  & 9772  &   518  & 950  & 262 &  227  & 179    \\ \hline
 458640299  & 12660  &   842  & 1539 &  423 &  366  & 289     \\ \hline
 698880299  & 13950  &   1087  & 1986 &  551  & 478  & 377    \\ \hline
 851760299  & 15904  &  1379 & 2524 &  697 &  604 &  475     \\ \hline
\end{tabular}
\end{center}
\caption{Memory requirements (in MB) for the storage of $M_{D,13}(x)$, Weber, $M_{D,5,7}(x)$, $M_{D,3, 13}(x)$ and Ramanujan polynomials.}
\label{memorytable}
\end{table}

Comparing the number of bits for the storage of all classes of polynomials, it is clear that the memory required
for the storage of the Ramanujan polynomials is smaller than the memory needed for the other three classes of polynomials.
The percentages are the same as in the precision requirements of the polynomials with one exception: Weber polynomials.
Notice that the degree of Weber polynomials 
is $3h$ and thus the memory used for the storage of Ramanujan polynomials is not only 42\% (like the precision
requirements) less than the
corresponding memory needed for the Weber polynomials but approximately 81\% less! This means that regarding the storage 
requirements of all polynomials, Weber polynomials are by far the worst choice. 
In Table~\ref{memorytable}
we present the memory in MB needed for the storage of all classes of polynomials for few values of $D$.
The difference in the efficiency of the construction of all classes of polynomials can be easily understood
noticing the polynomials for $D=299$ and $h=8$. Even though this is a small value for the discriminant,
the difference in the size of the coefficients of the polynomials is remarkable. 
In particular, 25 bits are required for the storage of the coefficients of the $T_{299}(x)$ polynomial,
188 bits for
the storage of $W_{299}(x)$ polynomial, 112 bits for $M_{299, 13}(x)$ polynomial, 
31 bits for $M_{299, 3, 13}(x)$ and 32 bits for $M_{299,5,7}(x)$.

\begin{equation*}
W_{299}(x) = x^{24} - 8x^{23}-12x^{22}-28x^{21}-56x^{20} -40x^{19} + 144x^{18} +144x^{17} +16x^{16} -112x^{15} -224x^{14} -416x^{13}
\end{equation*}
\begin{equation*}
-32x^{12} +256x^{11} +704x^{10} + 832x^{9} +640x^{8} -384x^{7} -1792x^{6} -1280x^{5} -256x^{4} +1280x^{3} +1536x^{2} +512x +256
\end{equation*}

\begin{equation*}
M_{299,13}(x) = x^{8} + 78x^{7}+793x^{6}+5070x^{5}+20956x^{4} +65910x^{3} + 134017x^{2} +171366x +28561
\end{equation*}

\begin{equation*}
M_{299,5,7}(x) = x^{8} - 8x^{7}+31x^{6}-22x^{5}+28x^{4} -2x^{3} - 19x^{2} +8x -1
\end{equation*}

\begin{equation*}
M_{299,3,13}(x) = x^{8} - 6x^{7}+16x^{6}+12x^{5}-23x^{4} +12x^{3} + 16x^{2} -6x +1
\end{equation*}

\begin{equation*}
T_{299}(x) = x^{8} + x^{7}-x^{6}-12x^{5}+16x^{4} -12x^{3} + 15x^{2} -13x +1
\end{equation*}

The time efficiency of the construction of the polynomials is clearly proportionate to the corresponding precision 
requirements. However, notice that computing the Weber and $M_{D,l}(x)$ polynomials amounts to $2h$ evaluations of 
the eta function $\eta$ while for Ramanujan and $M_{D, p_1, p_2}(x)$ polynomials we need to evaluate the function 
$3h$ and $4h$
times respectively. This could be a disadvantage for Ramanujan and $M_{D, p_1, p_2}(x)$ polynomials, but this is not the case.
In particular, it was shown in \cite{EM02} that is sufficient for any polynomial to precompute the values of $\eta$ only at the $h$ reduced quadratic forms. 
Finally, we note that the time required for the transformation of a root of a Weber, Ramanujan or $M_{D,l}(x)$
polynomial to a root of
the corresponding Hilbert polynomial is approximately the same. The situation gets worse when $M_{D, p_1, p_2}(x)$
polynomials are used, because the time for the transformation and the storage of the modular polynomials
are larger.

In conclusion, we showed that Ramanujan polynomials clearly outweight in every aspect all previously used class polynomials 
for all values of the
discriminant $D \equiv 3 \bmod 8$ and therefore their use is particularly favored in the CM method for the generation of
prime order ECs.


\begin{thebibliography}{99}

\bibitem{AM93} A.O.L. Atkin and F. Morain, Elliptic curves and
primality proving, {\em Mathematics of Computation}, 61(1993), pp.
29-67.

\bibitem{ACDFLNV06} R. M. Avanzi, H. Cohen, C. Doche, G. Frey, T. Lange, K. Nguyen, F. Vercauteren,
{\em Handbook of Elliptic and Hyperelliptic Curve Cryptography}, Chapman \& Hall/CRC, 2006.


\bibitem{B01} H.~Baier, Elliptic Curves of Prime Order over Optimal
Extension Fields for Use in Cryptography, in {\em Progress in
Cryptology} -- INDOCRYPT 2001, LNCS Vol.~2247, Springer-Verlag,
pp. 99-107, 2001.

\bibitem{B02} H.~Baier, Efficient Algorithms for Generating Elliptic
Curves over Finite Fields Suitable for Use in Cryptography,{\em  PhD Thesis,}
Dept.~of Computer Science, Technical Univ.~of Darmstadt, May 2002.

%\bibitem{BB00} H. Baier and J. Buchmann,
%Efficient construction of cryptographically strong elliptic curves,
%in {\em Progress in Cryptology} -- INDOCRYPT 2000,
%LNCS Vol.~1977, Springer-Verlag,
%pp. 191-202, 2000.

%\bibitem{B70} E. R. Berlekamp,
%Factoring polynomials over large finite fields, {\em Mathematics
%of Computation} 24(1970), pp. 713-735.

\bibitem{Berndt-Chan}
B.~C. Berndt and H.~H. Chan, \emph{Ramanujan and the modular
  {$j$}-invariant}, Canad. Math. Bull. \textbf{42} (1999), no.~4, 427--440.


\bibitem{BSS99} I. Blake, G. Seroussi, and N. Smart,
{\em Elliptic curves in cryptography} , London Mathematical
Society Lecture Note Series 265, Cambridge University Press, 1999.

\bibitem{BLS01} D. Boneh, B. Lynn, and H. Shacham, Short signatures from
the Weil pairing, in {\em ASIACRYPT 2001}, LNCS Vol.~2248,
Springer-Verlag, pp. 514-532, 2001.

\bibitem{BS06} R. Br\"{o}ker and P. Stevenhagen, Constructing elliptic curves of prime order, Preprint, Nov. 2006.

%\bibitem{BS07} R. Br\"{o}ker and P. Stevenhagen, Efficient CM-constructions of elliptic curves over finite fields, 
%{\em Mathematics of Computation}, \textbf{76}, 2007, pp. 2161--2179.


%\bibitem{GeeHuatTan}
%H.~H. Chan, A. Gee, and V. Tan, \emph{Cubic singular moduli,
%  {R}amanujan's class invariants {$\lambda\sb n$} and the explicit {S}himura
%  reciprocity law}, Pacific J. Math. \textbf{208} (2003), no.~1, 23--37.


\bibitem{C93} H. Cohen, { A Course in Computational Algebraic
Number Theory}, {\em Graduate Texts in Mathematics,} {\bf 138},
Springer-Verlag, Berlin, 1993.

\bibitem{C08} G. Cornacchia, Su di un
metodo per la risoluzione in numeri interi dell' equazione
$\sum_{h=0}^{n} C_{h}x^{n-h}y^h = P$, {\em Giornale di Matematiche
di Battaglini} 46 (1908), pp. 33-90.

\bibitem {C89} D. A. Cox,
{ Primes of the form $x\sp 2 + ny\sp 2$}, John Wiley and Sons, New York, 1989.

\bibitem{EM02}
A. Enge and F. Morain, Comparing invariants for  class fields of imaginary
quadratic fields, in {\em Algebraic Number Theory} -- ANTS V,
LNCS Vol.~2369, Springer-Verlag, pp. 252-266, 2002.

\bibitem{ES03}
A. Enge and R. Schertz, Constructing elliptic curves from modular curves of
positive genus, {\em Preprint,} 2003.

\bibitem{ES04}
A. Enge and R. Schertz, Constructing elliptic curves over finite fields using
double eta-quotients, {\em J. Th\'eor. Nombres Bordeaux}, \textbf{16} (2004), pp.555--568.

\bibitem{ES05}
A. Enge and R. Schertz, Modular curves of composite level,
{\em Acta Arithmetica}, 118 (2), (2005), pp.129--141.

\bibitem{FR94}
G. Frey and H.G. R\"{u}ck, A remark concerning $m$-divisibility and the discrete
logarithm problem in the divisor class group of curves,
{\em Mathematics of Computation}, 62(1994), pp.865--874.


%\bibitem{GMcKee00} S. Galbraith and J. McKee, The probability that the number of
%points on an elliptic curve over a finite field is prime, {\em
%Journal of the London Mathematical Society}, 62(2000), no. 3, pp.
%671-684.

\bibitem{GeeBordeaux}
A. Gee, {Class invariants by {S}himura's reciprocity law}, {\em J. Th\'eor.
  Nombres Bordeaux} \textbf{11} (1999), no.~1, 45--72, Les XX\`emes Journ\'ees
  Arithm\'etiques (Limoges, 1997).

\bibitem{GeeStevenhagen}
A. Gee and P. Stevenhagen, {Generating class fields using {S}himura
  reciprocity}, {\em Algorithmic number theory (Portland, OR, 1998)}, LNCS Vol.~1423, 
Springer-Verlag, pp.~441--453, 1998.


\bibitem{gnu:gnu} GNU multiple precision library, edition 4.2.1,
 2007. Available at: {\tt http://www.swox.com/gmp}.

\bibitem{HinSil}
M. Hindry, and  J. Silverman,
 {\em Diophantine geometry An introduction}, Graduate Texts in Mathematics, Springer-Verlag, New York, 2000. 



\bibitem{ieee} IEEE P1363/D13,
{\em Standard Specifications for Public-Key Cryptography}, 1999.
{\tt http://grouper.ieee.org/groups/1363/tradPK/draft.html}.

%\bibitem{KY89}
%E. Kaltofen and N. Yui, Explicit construction of the Hilbert class
%fields of imaginary quadratic fields by integer lattice reduction.
%Research Report 89-13, Rensselaer Polytechnic Institute, May 1989.

\bibitem{KVY89}
E. Kaltofen, T. Valente, and N. Yui, An Improved Las Vegas
Primality Test, in {\em Proc.~ACM-SIGSAM 1989
International Symposium on Symbolic and Algebraic Computation},
pp. 26-33, 1989.

%\bibitem{KSZ02_esa}
%E.~Konstantinou, Y.~Stamatiou, and C.~Zaroliagis,
%A Software Library for Elliptic Curve Cryptography,
%in {\em Proc.~10th European Symposium on Algorithms}
%-- ESA 2002 (Engineering and Applications Track),
%LNCS Vol.~2461, Springer-Verlag, pp. 625-637, 2002.

\bibitem{KSZ04_icisc}
E.~Konstantinou, A.~Kontogeorgis, Y.~Stamatiou, and C.~Zaroliagis, Generating Prime Order Elliptic Curves:
Difficulties and Efficiency Considerations, in {\em
International Conference on Information Security and Cryptology} -- ICISC 2004,
LNCS Vol.~3506, Springer-Verlag, pp. 261-278, 2005.

\bibitem{KonstKonto}
E.~Konstantinou, A.~Kontogeorgis, Computing Polynomials of the Ramanujan $t_n$ Class Invariants,
{\em Arxiv:math.NT/0610372} v1, 11 Oct. 2006, to appear in Canadian Mathematical Bulletin.


\bibitem{KSZ02}
E.~Konstantinou, Y.~Stamatiou, and C.~Zaroliagis,
On the Efficient Generation of Elliptic Curves over Prime Fields,
in {\em Cryptographic Hardware and Embedded Systems}
-- CHES 2002, LNCS Vol.~2523,
Springer-Verlag, pp. 333-348, 2002.

%\bibitem{KSZ03} E.~Konstantinou, Y.C.~Stamatiou, and C.~Zaroliagis, On the
%Construction of Prime Order Elliptic Curves, in {\em
%Progress in Cryptology} -- INDOCRYPT 2003,
%LNCS Vol.~2904, Springer-Verlag, pp. 309-322, 2003.


\bibitem{LZ94}
G.J.~Lay and H.~Zimmer, Constructing Elliptic Curves with Given
Group Order over Large Finite Fields, in {\em Algorithmic Number
Theory} -- ANTS-I, LNCS Vol.~877,
Springer-Verlag, pp. 250-263, 1994.

%\bibitem{lidia} LiDIA.
%{\em A library for computational number theory},
%Technical University of Darmstadt. Available from
%{\tt http://www.informatik.tu-darmstadt.de/TI/LiDIA/Welcome.html}.

\bibitem{MOV93} A. J. Menezes, T. Okamoto and S. A. Vanstone,
Reducing elliptic curve logarithms to a finite field, {\em IEEE
Trans. Info. Theory}, 39(1993), pp. 1639-1646.

\bibitem{M92}
A. Miyaji, Elliptic curves over $F_p$ suitable for cryptosystems, in {\em
Advances in Cryptology} -- AUSCRYPT 1992,
LNCS Vol.~718, Springer-Verlag, pp. 492-504, 1992.


\bibitem{MNT00}
A. Miyaji, M. Nakabayashi, and S. Takano, Characterization of Elliptic Curve Traces
under FR-reduction, in {\em
International Conference on Information Security and Cryptology} -- ICISC 2000,
LNCS Vol.~2015, Springer-Verlag, pp. 90-108, 2001.

\bibitem{MNT01}
A. Miyaji, M. Nakabayashi, and S. Takano, New explicit conditions of elliptic curve
traces for FR-reduction, {\em IEICE Transactions on Fundamentals}, E84-A(5):1234-1243,
2001.

\bibitem{M90} F.~Morain, Construction of hilbert class fields of imaginary quadratic fields and
dihedral equation modulo $p$, Report 1087, INRIA, 1989.


\bibitem{M00}
F.~Morain, Modular curves and class invariants, {\em Preprint,} June 2000.

%\bibitem{M02}
%F.~Morain, Computing the cardinality of CM elliptic curves using
%torsion points, Preprint, October 2002.

\bibitem{MP97}
V. M\"{u}ller and S. Paulus, On the Generation of
Cryptographically Strong Elliptic Curves, preprint, 1997.


\bibitem{NM05}
Y.~Nogami and Y.~Morikawa, A method for distinguishing the two candidate elliptic curves
in CM method, in {\em
International Conference on Information Security and Cryptology} -- ICISC 2004,
LNCS Vol.~3506, Springer-Verlag, pp. 249-260, 2005.


%\bibitem{NM03} Y. Nogami and Y. Morikawa, Fast generation of elliptic curves
%with prime order over $F_{p^{2^c}}$, in {\em Proc. of the International
%workshop on Coding and Cryptography}, March 2003.

\bibitem{PARI2} PARI/GP, version {\tt 2.3.1}, Bordeaux, 2005.
Available at: {\tt http://pari.math.u-bordeaux.fr/}. 

\bibitem{PH78} G. C. Pohlig and M. E. Hellman,
An improved algorithm for computing logarithms over $GF(p)$ and
its cryptographic significance, {\em IEEE Trans. Info. Theory}, 24
(1978), pp. 106-110.

\bibitem{RamNotebooks}
S. Ramanujan, {Notebooks. {V}ols. 1, 2}, {\em Tata Institute of
  Fundamental Research}, Bombay, 1957.

\bibitem{RS07} K. Rubin and A. Silverberg, Choosing the correct elliptic curve in the CM method, 
preprint, 2007.

\bibitem{Sarnak}
P. Sarnak, {\em  Selberg's eigenvalue conjecture}. 
 Notices Amer. Math. Soc. 42 (1995), no. 11, 1272--1277. 


\bibitem{SA98} T. Satoh and K. Araki,
Fermat quotients and the polynomial time discrete log algorithm
for anomalous elliptic curves, {\em Comm. Math. Univ. Sancti
Pauli}, 47(1998), pp. 81-91.

\bibitem{SSK01} E.~Sava\c{s}, T.A.~Schmidt, and \c{C}.K.~Ko\c{c},
Generating Elliptic Curves of Prime Order, in {\em Cryptographic
Hardware and Embedded Systems} -- CHES 2001, LNCS Vol.~2162,
Springer-Verlag, pp. 145-161, 2001.

\bibitem{S02}
R. Schertz, Weber's class invariants revisited, {\em Journal de Th\'{e}orie des Nombres de Bordeaux},
4(2002), pp. 325-343, 2002.

\bibitem{S95} R. Schoof,
Counting points on elliptic curves over finite fields,
{\em J.~Theorie des Nombres de Bordeaux}, 7(1995), pp. 219-254, 1995.

\bibitem{SB04}
M. Scott and P. S.L.M. Barreto, Generating more MNT elliptic curves, {\em Cryptology ePrint Archive},
Report 2004/058, 2004.

\bibitem{S86} J. H. Silverman, {\em The Arithmetic of Elliptic Curves},
Springer-Verlag, GTM 106, 1986.

\bibitem{S04} I. Stewart, {\em Galois Theory}, Third Edition, Chapman \& Hall/CRC,
Boca Raton, FL, 2004.

\bibitem{ST87} I. Stewart and D. Tall, {\em Algebraic Number Theory}, Second Edition,
Chapman \& Hall, London, 1987.

%\bibitem{V92} T.~Valente,
%{\em A distributed approach to proving large numbers prime},
%Rensselaer Polytechnic Institute Troy, New York, PhD Thesis, August 1992.

\bibitem{WE01} A. Weng, {\em Konstruktion kryptographisch geeigneter Kurven
mit komplexer Multiplikation}, PhD thesis, Institut f\"{u}r
Experimentelle Mathematik, Universit\"{a}t GH Essen, 2001.



\end{thebibliography}
\end{document}